\documentclass[a4paper]{amsart}

\usepackage[T1]{fontenc}
\usepackage[latin1]{inputenc}
\usepackage{amssymb}
\usepackage{euscript}
\usepackage{amsxtra}
\usepackage{ae}
\usepackage[all]{xy}
\usepackage{enumerate}
\usepackage{amstext}
\usepackage{mathrsfs}

\usepackage[lite]{amsrefs}

\theoremstyle{plain}
\newtheorem{thm}{Theorem}
\newtheorem{lemma}[thm]{Lemma}
\newtheorem{prop}[thm]{Proposition}
\newtheorem{cor}[thm]{Corollary}

\theoremstyle{definition}
\newtheorem{defn}[thm]{Definition}

\theoremstyle{remark}
\newtheorem{rmk}[thm]{Remark}
\newtheorem{note}[thm]{Note}
\newtheorem{ex}[thm]{Example}

\newcommand*{\ebar}{{\underline{E}G}}
\newcommand*{\C}{{\mathbb C}}
\newcommand*{\Z}{{\mathbb Z}}
\newcommand*{\Ztwo}{{\mathbb Z/2}}
\newcommand*{\R}{{\mathbb R}}
\newcommand*{\N}{{\mathbb N}}
\newcommand*{\Q}{{\mathbb Q}}

\newcommand*{\X}{\mathscr{X}}
\newcommand*{\Y}{\mathscr{Y}}
\newcommand*{\rips}{\mathscr{P}}

\newcommand*{\curD}{{\mathfrak{D}}}
\newcommand*{\relD}{{\mathfrak{D}_{\mathrm{red}}}}

\newcommand*{\vv}{\mathfrak{\bar{c}}}
\newcommand*{\vvr}{\mathfrak{\bar{c}^{red}} }
\newcommand*{\bv}{\mathfrak{c}}
\newcommand*{\bvr}{\mathfrak{c^{red}}}
\newcommand*{\Borelvv}{\overline{\mathfrak{b}}}
\newcommand*{\Borelbv}{\mathfrak{b}}

\newcommand*{\Proj}{{\mathsf{P}}}
\newcommand*{\Dirac}{{\mathsf{D}}}

\newcommand*{\KK}{\mathrm{KK}}

\newcommand*{\K}{\mathrm{K}}

\newcommand*{\KX}{\mathrm{KX}}

\newcommand*{\Var}{\mathrm{Var}}
\newcommand*{\ev}{\mathrm{ev}}

\DeclareMathOperator{\Coker}{Coker}

\DeclareMathOperator{\supp}{supp}

\newcommand*{\ripsn}{\mathrm{P}}

\newcommand*{\Mult}{{\mathcal{M}}}
\newcommand*{\Calkin}{{\mathcal{Q}}}

\newcommand*{\Comp}{{\mathbb{K}}}
\newcommand*{\Bound}{{\mathbb{B}}}

\newcommand*{\brd}{-\hspace{0pt}}
\newcommand*{\nbd}{\nobreakdash-\hspace{0pt}}

\newcommand*{\abs}[1]{\lvert#1\rvert}

\newcommand*{\braket}[2]{\langle#1,#2\rangle}
\newcommand*{\norm}[1]{\lVert#1\rVert}

\newcommand*{\cross}{\mathbin{\ltimes}}

\newcommand*{\defeq}{\mathrel{:=}}

\begin{document}

\title{Dualizing the coarse assembly map}
\author{Heath Emerson}
\email{hemerson@math.uni-muenster.de}

\author{Ralf Meyer}
\email{rameyer@math.uni-muenster.de}

\address{Mathematisches Institut\\
         Westf.\ Wilhelms-Universität Münster\\
         Ein\-stein\-stra\-ße 62\\
         48149 Münster\\
         Germany}

\begin{abstract}
  We formulate and study a new coarse (co-)assembly map.  It involves
  a modification of the Higson corona construction and produces a map
  dual in an appropriate sense to the standard coarse assembly map.
  The new assembly map is shown to be an isomorphism in many cases.
  For the underlying metric space of a group, the coarse co-assembly
  map is closely related to the existence of a dual Dirac morphism and
  thus to the Dirac dual Dirac method of attacking the Novikov
  conjecture.
\end{abstract}


\thanks{This research was supported by the EU-Network \emph{Quantum
  Spaces and Noncommutative Geometry} (Contract HPRN-CT-2002-00280)
  and the \emph{Deutsche Forschungsgemeinschaft} (SFB 478).}

\maketitle

\section{Introduction}
\label{sec:intro}

It is shown in~\cite{EmersonMeyer} that a torsion free discrete
group~$G$ with compact classifying space~$BG$ has a dual Dirac
morphism (in the sense of~\cite{MeyerNest}) if and only if a certain
coarse co\brd{}assembly map
$$
\mu^*\colon \K_{*+1}(\bv(G)) \to \KX^*(G)
$$
is an isomorphism.  The $C^*$\nbd{}algebra $\bv(G)$ is called the
stable Higson corona of~$G$ and up to isomorphism only depends on the
coarse structure of~$G$.  The $\Ztwo$\nbd{}graded Abelian group
$\KX^*(G)$ is called the coarse $\K$\nbd{}theory of~$G$ and also
depend only on the coarse structure of~$G$.  Essentially the same
result holds for torsion free, countable, discrete groups with finite
dimensional~$BG$.  In particular, for this class of groups, the
existence of a dual Dirac morphism is a geometric invariant of~$G$.
In this article we introduce and study the map~$\mu^*$ in detail.  In
particular, we
\begin{enumerate}[(1)]
\item examine the relationship between~$\mu^*$ and the ordinary coarse
  Baum-Connes assembly map, and between $\bv(G)$ and compactifications
  of~$G$;
  
\item establish isomorphism of~$\mu^*$ for scalable spaces; 
  
\item establish isomorphism of~$\mu^*$ for groups which uniformly
  embed in Hilbert space.

\end{enumerate}

The stable Higson corona $\bv(X)$ has better functoriality properties
than the $C^*$\nbd{}algebra $C^*(X)$ that figures in the usual coarse
Baum-Connes assembly map (see
\cites{HigsonRoe,Roe,SkandalisYuTu,Yu2,Yu:BC_coarse_BC,Yu,Yu:embeddable}).
The assignment $X\mapsto\bv(X)$ is functorial from the coarse category
of coarse spaces to the category of $C^*$\nbd{}algebras and
$C^*$\nbd{}algebra homomorphisms.  The analogous statement for the
coarse $C^*$\nbd{}algebra $C^*(X)$ is true only after passing to
$\K$\nbd{}theory.  Moreover, the $C^*$\nbd{}algebra $\bv(X)$ is
designed to be closely related to certain bivariant Kasparov groups.
This is the source of a homotopy invariance result, which implies our
assertion for scalable spaces.  Another advantage of the stable Higson
corona and the coarse co\brd{}assembly map is their relationship with
alternative approaches to the Novikov conjecture, namely, almost flat
$\K$\nbd{}theory (see \cites{CGM1, EmersonKaminker}) and the Lipschitz
approach of~\cite{CGM2}.

The map~$\mu^*$ is an isomorphism for any discrete group~$G$ that has
a dual Dirac morphism, without any hypothesis on~$BG$.  This is how we
are going to prove isomorphism of~$\mu^*$ for groups that uniformly
embed in a Hilbert space: we show that such groups have a dual Dirac
morphism.  Actually, already the existence of an approximate dual
Dirac morphism implies that~$\mu^*$ is an isomorphism.  Using results
of Gennadi Kasparov and Georges Skandalis
(\cite{Kasparov-Skandalis:Bolic}), it follows that the coarse
co-assembly map is an isomorphism for groups acting properly by
isometries on bolic spaces.  The usual coarse Baum-Connes conjecture
for a group~$G$ is equivalent to the Baum-Connes conjecture with
coefficients $\ell^\infty(G)$ (\cite{Yu:BC_coarse_BC}).  Despite this,
it is not known whether the existence of an action of~$G$ on a bolic
space implies the coarse Baum-Connes conjecture for~$G$.  The
existence of a dual Dirac morphism only implies split injectivity of
the coarse Baum-Connes assembly map.

Given the above observations, we expect the coarse co-assembly to
become a useful tool in connection with the Novikov conjecture.
However, at the moment we have no examples of groups for which our
method proves the Novikov conjecture while others fail.  We also
remark that we do not know whether the map~$\mu^*$ is an isomorphism
for the standard counter-examples to the coarse Baum-Connes
conjecture.

Finally, we would like to thank the referee for his thorough report
and his useful comments.

\section{Coarse spaces}
\label{sec:coarse}

We begin by recalling the notion of a coarse space and some related
terminology (see~\cites{HigsonRoe, SkandalisYuTu}).  Then we introduce
$\sigma$\nbd{}coarse spaces, which are useful to deal with the Rips
complex construction.

Let~$X$ be a set.  We define the \emph{diagonal}~$\Delta_X$, the
\emph{transpose} of $E\subseteq X\times X$, and the \emph{composition}
of $E_1,E_2\subseteq X\times X$ by
\begin{align*}
  \Delta_X &\defeq \{(x,x)\in X\times X \mid x\in X\},
  \\
  E^t &\defeq \{(y,x)\in X\times X \mid (x,y)\in E\},
  \\
  E_1 \circ E_2 &\defeq
  \{(x,z) \in X\times X \mid
    \text{$(x,y)\in E_1$ and $(y,z)\in E_2$ for some $y\in X$}
  \}.
\end{align*}

\begin{defn}  \label{definitionofcoarsestructure}
  A \emph{coarse structure} on~$X$ is a collection~$\mathcal{E}$ of
  subsets $E\subseteq X\times X$---called \emph{entourages} or
  \emph{controlled subsets}---which satisfy the following axioms:
  \begin{enumerate}[\ref{definitionofcoarsestructure}.1.]
  \item if $E\in\mathcal{E}$ and $E'\subseteq E$, then
    $E'\in\mathcal{E}$ as well;

  \item if $E_1,E_2\in\mathcal{E}$, then $E_1\cup E_2\in\mathcal{E}$;

  \item if $E\in\mathcal{E}$ then $E^t\in\mathcal{E}$;

  \item if $E_1,E_2\in\mathcal{E}$, then $E_1\circ E_2\in\mathcal{E}$;

  \item $\Delta_X\in\mathcal{E}$;

  \item all finite subsets of $X\times X$ belong to~$\mathcal{E}$.

  \end{enumerate}
  
  A subset~$B$ of~$X$ is called \emph{bounded} if $B\times B$ is an
  entourage.  A collection of bounded subsets $(B_i)$ of~$X$ is called
  \emph{uniformly bounded} if $\bigcup B_i \times B_i$ is an
  entourage.

  A topology and a coarse structure on~$X$ are called
  \emph{compatible} if
  \begin{enumerate}[\ref{definitionofcoarsestructure}.1.]
  \setcounter{enumi}{6}
  \item some neighborhood of $\Delta_X\subseteq X\times X$ is an
    entourage;

  \item every bounded subset of~$X$ is relatively compact.

  \end{enumerate}
  
  A \emph{coarse space} is a locally compact topological space
  equipped with a compatible coarse structure.
\end{defn}

Since the intersection of a family of coarse structures is again a
coarse structure, we can define the \emph{coarse structure generated
  by} any set of subsets of $X\times X$.  We call a coarse structure
\emph{countably generated} if there is an increasing sequence of
entourages $(E_n)$ such that any entourage is contained in~$E_n$ for
some $n\in\N$.  If~$X$ is a coarse space and $Y\subseteq X$ is a
closed subspace, then $\mathcal{E}\cap (Y\times Y)$ is a coarse
structure on~$Y$ called the \emph{subspace coarse structure}.

Let~$X$ be a coarse space.  Then the closure of an entourage is again
an entourage.  Hence the coarse structure is already generated by the
closed entourages.  Using an entourage that is a neighborhood of the
diagonal, we can construct a uniformly bounded open cover of~$X$.
This open cover has a subordinate partition of unity because~$X$ is
locally compact.  We shall frequently use this fact.

\begin{ex}  \label{themetriccoarsestructure}
  Let $(X,d)$ be a metric space.  The \emph{metric coarse structure
  on~$X$} is the countably generated coarse structure generated by the
  increasing sequence of entourages
  $$
  E_R \defeq \{(x,y)\in X\times X \mid d(x,y) \le R\},
  \qquad R\in\N.
  $$
  A subset $E\subseteq X\times X$ is an entourage if and only if
  $d\colon X\times X\to\R_+$ is bounded on~$E$.  Note that this coarse
  structure depends only on the quasi-isometry class of~$d$.
  The metric~$d$ also defines a topology on~$X$.  This topology and
  the coarse structure are compatible and thus define a coarse space
  if and only if bounded subsets of~$X$ are relatively compact.  We
  call $(X,d)$ a \emph{coarse metric space} if this is the case.

  Conversely, one can show that any countably generated coarse
  structure on a set~$X$ can be obtained from some metric on~$X$ as
  above.  However, if~$X$ also carries a topology, it is not clear
  whether one can find a metric that generates both the coarse
  structure and the topology.
\end{ex}

\begin{ex}  \label{ex:group_coarse_structure}
  Any locally compact group~$G$ has a canonical coarse structure that
  is invariant under left translations.  It is generated by the
  entourages
  $$
  E_K \defeq \{(g_1,g_2)\in G\times G \mid g_1^{-1}g_2 \in K\},
  $$
  where~$K$ runs through the compact subsets of~$G$.  Together with
  the given locally compact topology on~$G$, this turns~$G$ into a
  coarse space.
\end{ex}

A \emph{coarse map} $\phi\colon X \to Y$ between coarse spaces is a
Borel map which maps entourages in~$X$ to entourages in~$Y$ and which
is \emph{proper} in the sense that inverse images of bounded sets are
bounded.  Two coarse maps $\phi,\psi\colon X\to Y$ are called
\emph{close} if $(\phi\times\psi)(\Delta_X)\subseteq Y\times Y$ is an
entourage.  The \emph{coarse category of coarse spaces} is the
category whose objects are the coarse spaces and whose morphisms are
the equivalence classes of coarse maps, where we identify two maps if
they are close.  A coarse map is called a \emph{coarse equivalence} if
it is an isomorphism in this category.  Two coarse spaces are called
\emph{coarsely equivalent} if they are isomorphic in this category.

\begin{lemma}  \label{lem:discretization}
  Let~$X$ be a countably generated coarse space.  Then there exists a
  countable discrete subset $Z\subseteq X$ such that the inclusion
  $Z\to X$ is a coarse equivalence.  Here we equip~$Z$ with the
  subspace coarse structure and the discrete topology.  Thus~$X$ is
  coarsely equivalent to a countably generated, discrete coarse space.
\end{lemma}

\begin{proof}
  We claim that any countably generated coarse space is
  $\sigma$\nbd{}compact.  To see this, fix a point $x_0\in X$ and an
  increasing sequence~$(E_n)$ of entourages that defines the coarse
  structure.  The sets $K_n\defeq \{x\in X\mid (x,x_0)\in E_n\}$ are
  bounded and hence relatively compact.  Their union is all of~$X$
  because $\bigcup E_n= X\times X$.  Thus~$X$ is
  $\sigma$\nbd{}compact.  Let~$(B_i)$ be a uniformly bounded open
  cover of~$X$.  By $\sigma$\nbd{}compactness, we can choose a
  countable partition of unity $(\rho_n)$ on~$X$ subordinate to this
  covering.  Let $B'_n\defeq\rho_n^{-1}((0,\infty))$.  These sets form
  a countable, locally finite, uniformly bounded open covering of~$X$.
  Choose a point~$x_n$ in each~$B'_n$.  The subset $Z\defeq \{x_n\}$
  has the required properties.
\end{proof}

We shall also use formal direct unions of coarse spaces, which we call
$\sigma$\nbd{}coarse spaces.  Let $(X_n)_{n\in\N}$ be an increasing
sequence of subsets of a set~$\X$ with $\X=\bigcup X_n$ such that
each~$X_n$ is a coarse space and the coarse structure and topology
on~$X_m$ are the subspace coarse structure and topology from~$X_n$ for
any $n\ge m$.  Then we can equip~$\X$ with the direct limit topology
and with the coarse structure that is generated by the coarse
structures of the subspaces~$X_n$.  This coarse structure is
compatible with the topology, but the topology need not be locally
compact.  In this situation, we call the inductive system of coarse
spaces~$(X_n)$ or, by abuse of notation, its direct limit~$\X$ a
\emph{$\sigma$\nbd{}coarse space}.  We define $\sigma$\nbd{}locally
compact spaces similarly.  Notice that the system $(X_n)$ is part of
the structure of~$\X$ even if the direct limit topology on~$\X$ is
locally compact.

\begin{ex}  \label{ex:metric_rips}
  Let $(X,d)$ be a discrete metric space with the property that
  bounded subsets are finite.  Let $\ripsn_n(X)$ denote the set of
  probability measures on~$X$ whose support has diameter at most~$n$.
  This is a locally finite simplicial complex and hence a locally
  compact topological space.  We equip $\ripsn_n(X)$ with the coarse
  structure generated by the increasing sequence of entourages
  $$
  \{(\mu,\nu)\in\ripsn_n(X)\times\ripsn_n(X)\mid
  \supp\mu\times\supp\nu\subseteq E_R\}
  $$
  for $R\in\N$, with~$E_R$ as in
  Example~\ref{themetriccoarsestructure}.  This turns $\ripsn_n(X)$
  into a coarse space for any $n\in\N$ and turns $\rips_X \defeq
  \bigcup \ripsn_n (X)$ into a $\sigma$\nbd{}coarse space.  We discuss
  this example in greater detail in
  Section~\ref{sec:coarse_coassembly}.
\end{ex}

\begin{ex}  \label{ex:proper_action_coarse}
  Let~$G$ be a second countable, locally compact group and let~$X$ be
  a \emph{$G$\nbd{}compact}, proper $G$\nbd{}space.  We equip~$X$ with
  the coarse structure that is generated by the $G$\nbd{}invariant
  entourages
  $$
  E_L \defeq \bigcup_{g\in G} gL\times gL,
  \qquad
  \text{$L\subseteq X$ compact.}
  $$
  Since~$X$ is necessarily $\sigma$\nbd{}compact, this coarse
  structure is countably generated.  It is also compatible with the
  topology of~$X$, so that~$X$ becomes a countably generated coarse
  space.  For $X=G$ with the action by left translation, this
  reproduces the coarse structure of
  Example~\ref{ex:group_coarse_structure}.  For any $x\in X$, the
  orbit map $G\to X$, $g\mapsto g\cdot x$, is a coarse equivalence,
  and these maps for different points in~$X$ are close.  Thus we
  obtain a canonical isomorphism $X\cong G$ in the coarse category of
  coarse spaces.
  
  Let $\ebar$ denote a second countable, not necessarily locally
  compact model for the classifying space for proper actions of~$G$ as
  in~\cite{Baum-Connes-Higson}.  We can write $\ebar$ as an increasing
  union of a sequence of $G$\nbd{}compact, $G$\nbd{}invariant, closed
  subsets $X_n\subseteq\ebar$.  Turning each~$X_n$ into a coarse space
  as above, we turn $\ebar$ into a $\sigma$\nbd{}coarse space.
\end{ex}

In the above two examples, the maps $X_n\to X_{n+1}$ are coarse
equivalences.  This happens in all examples that we need, and the
general case is more difficult.  Therefore, \emph{we restrict
attention in the following to $\sigma$\nbd{}coarse spaces for which
the maps $X_n\to X_{n+1}$ are coarse equivalences.}

\section{Functions of vanishing variation}
\label{sec:vanishing_variation}

Let~$X$ be a coarse space and let~$D$ be a $C^*$\nbd{}algebra.  We
define the $C^*$\nbd{}algebras $\vv(X,D)$ and $\bv(X,D)$ and discuss
their relationship to the Higson compactification and the Higson
corona and to admissible compactifications.  Then we investigate their
functoriality properties.

\begin{defn}  \label{definitionofvanishingvariation}
  Let~$\X$ be a $\sigma$\nbd{}coarse space and let~$Y$ be
  a metric space.  Let $f\colon \X \to Y$ be a Borel map (that is,
  $f|_{X_n}$ is a Borel map for all $n\in\N$).  For an entourage
  $E\subseteq X_n\times X_n$, $n\in\N$, we define
  $$
  \Var_E\colon X_n\to [0,\infty),
  \qquad
  \Var_E f(x) \defeq
  \sup \{ d\bigl(f(x),f(y)\bigr) \mid (x,y)\in E \}.
  $$
  We say that~$f$ has \emph{vanishing variation} if $\Var_E$ vanishes
  at~$\infty$ for any such~$E$, that is, for any $\epsilon>0$ the set
  of $x\in X_n$ with $\Var_E f(x)\ge\epsilon$ is bounded.
\end{defn}

If the coarse structure comes from a metric~$d$ on~$X$, we also let
$$
\Var_R f(x) \defeq \sup \{d(f(x),f(y)) \mid d(x,y)\le R\}
$$
for $R\in\R_+$.  This is the variation function associated to the
entourage~$E_R$ defined in Example~\ref{themetriccoarsestructure}.
Hence we can also use the functions $\Var_R f$ to define vanishing
variation.

\begin{defn}  \label{definitionofthestablehigsoncorona}
  For any coarse space~$X$ and any $C^*$\nbd{}algebra $D$, we let
  $\vv(X,D)$ be the $C^*$\nbd{}algebra of bounded, continuous
  functions of vanishing variation $X\to D\otimes\Comp$.  Here~$\Comp$
  denotes the $C^*$\nbd{}algebra of compact operators on a separable
  Hilbert space.  The quotient $\bv(X,D) \defeq \vv(X,D)/
  C_0(X,D\otimes\Comp)$ is called the \emph{stable Higson corona
    of~$X$ with coefficients~$D$.}
\end{defn}

\begin{note}
  When $D=\C$ we abbreviate $\vv(X,D)$ and $\bv(X,D)$ to $\vv(X)$ and
  $\bv(X)$, respectively, and call $\bv(X)$ the \emph{stable Higson
  corona of~$X$}.
\end{note}

The reason for our terminology is the analogy with the Higson corona
constructed in~\cite{Hig3}.  Let~$X$ be a coarse metric space.  The
\emph{Higson compactification}~$\eta X$ of~$X$ is the maximal ideal
space of the $C^*$\nbd{}algebra of continuous, bounded functions
$X\to\C$ of vanishing variation.  The \emph{Higson corona} of~$X$ is
$\partial_\eta X\defeq \eta X \setminus X$.  By construction,
$M_n(\C)\otimes C(\eta X)=C(\eta X,M_n)$ is the $C^*$\nbd{}algebra of
bounded, continuous functions $X\to M_n(\C)$ of vanishing variation
and $C(\partial_\eta X,M_n(\C)) = C(\eta X,M_n(\C))/C_0(X,M_n(\C))$.
Of course, these $C^*$\nbd{}algebras are contained in $\vv(X)$ and
$\bv(X)$, respectively.  Since $\Comp=\varinjlim M_n(\C)$, we also
obtain canonical embeddings
$$
\Comp\otimes C(\eta X)\cong C(\eta X,\Comp) \subseteq \vv(X),
\qquad
\Comp\otimes C(\partial_\eta X)\cong C(\partial_\eta X,\Comp)
\subseteq \bv(X).
$$
Similarly, we obtain embeddings
$$
C(\eta X,D\otimes \Comp) \subseteq \vv(X,D),
\qquad
C(\partial_\eta X,D\otimes \Comp)
\subseteq \bv(X,D)
$$
for any $C^*$\nbd{}algebra~$D$.  It turns out that $\vv(X)$ is
strictly larger than $C(\eta X,\Comp)$.  If $f\in C(\eta X,\Comp)$,
then $f(X)\subseteq f(\eta X)$ must be a relatively compact subset
of~$\Comp$.  Conversely, one can show that a continuous function
$X\to\Comp$ of vanishing variation with relatively compact range
belongs to $C(\eta X,\Comp)$.  However, functions in $\vv(X)$ need not
have relatively compact range.

It is often preferable to replace the Higson compactification by
smaller ones that are metrizable.  This is the purpose of the
following definition.

\begin{defn}[\cite{Hig2}]
  Let~$X$ be a metric space and let $i\colon X\to Z$ be a metrizable
  compactification of~$X$.  We call~$Z$ \emph{admissible} if there
  is a metric on~$Z$ generating the topology on~$Z$ for which the
  inclusion $i\colon X\to Z$ has vanishing variation.
\end{defn}

\begin{ex}
  The following are examples of admissible compactifications:
  \begin{enumerate}[(1)]
  \item the one-point compactification of an arbitrary metric space;
    
  \item the hyperbolic compactification of a metric space that is
    hyperbolic in the sense of Gromov; 
    
  \item the visibility compactification of a complete, simply
    connected, non\brd{}positively curved manifold.

  \end{enumerate}
\end{ex}

\begin{prop}  \label{pro:admissible_compactification}
  Let~$X$ be a metric space and let $i\colon X\to Z$ be an admissible
  compactification of~$X$.  Then there are canonical injective
  $*$\nbd{}homomorphisms
  \begin{gather*}
    C(Z,D\otimes\Comp) \to
    C(\eta X,D\otimes\Comp) \to \vv(X,D),
    \\
    C(Z\setminus X,D\otimes\Comp) \to
    C(\partial_\eta X,D\otimes\Comp) \to \bv(X,D)
  \end{gather*}
  for any $C^*$\nbd{}algebra~$D$.  Any class in
  $\K_*(C(\partial_\eta X,D\otimes\Comp))$ is the image of a class
  in $\K_*(C(Z\setminus X,D\otimes\Comp))$ for some admissible
  compactification~$Z$.
\end{prop}

\begin{proof}
  The admissible compactifications of~$X$ are exactly the metrizable
  quotients of the Higson compactification~$\eta X$.  Thus the
  $C^*$\nbd{}algebras $C(Z)$ for admissible compactifications of~$X$
  are exactly the separable $C^*$\nbd{}subalgebras of $C(\eta X)$.
  This implies the corresponding assertion for $C(Z\setminus X)$ and
  $C(\partial_\eta X)$ and also for tensor products with
  $D\otimes\Comp$.  The assertion about $\K$\nbd{}theory follows
  because any $C^*$\nbd{}algebra is the inductive limit of its
  separable $C^*$\nbd{}subalgebras and $\K$\nbd{}theory commutes with
  inductive limits.
\end{proof}

Even on the level of $\K$\nbd{}theory, $C(\eta X,\Comp)$ and $\vv(X)$
are drastically different.  The induced map $\K_*(C(\eta X,\Comp)) \to
\K_*(\vv(X))$ is uncountable-to-one already in rather simple examples,
as we shall see in Section~\ref{sec:vanishing}.  The map $\K_*(C(\eta
X,\Comp)) \to \K_*(\vv(X))$ may also fail to be surjective.  For
instance, this happens for the well-spaced ray (see
Example~\ref{dualcoarsebaumconnesforwellspacedray} below).  However,
we do not know of a uniformly contractible example with this property.
If the map $\K_*(C(\eta X,\Comp))\to\K_*(\vv(X))$ is surjective, then
Proposition~\ref{pro:admissible_compactification} yields that any
class in $\K_*(\vv(X))$ can already be realized on some admissible
compactification.

Now we turn to the functoriality of the algebras $\vv(X,D)$ and
$\bv(X,D)$ with respect to the coarse space~$X$.  The functoriality in
the coefficient algebra~$D$ is analyzed in
Section~\ref{sec:functor_in_coefficient}.  It is evident that
$\vv(X,D)$ and $\bv(X,D)$ and the extension
$$
0 \to C_0(X,D\otimes\Comp)\to\vv(X,D)\to\bv(X,D)\to0
$$
are functorial for \emph{continuous} coarse maps $X\to X'$.  Of
course, the morphisms in the category of $C^*$\nbd{}algebras are the
$*$\nbd{}homomorphisms.  We can drop the continuity hypothesis for
$\bv(X,D)$:

\begin{prop}  \label{functorialityofcoronaalgebra}
  Let~$D$ be a $C^*$\nbd{}algebra and let $X$ and~$Y$ be coarse
  spaces.  A coarse map $\phi\colon X\to Y$ induces a
  $*$\nbd{}homomorphism $\phi^*\colon \bv(Y,D)\to\bv(X,D)$.  Close
  maps induce the same $*$\nbd{}homomorphism $\bv(Y,D)\to\bv(X,D)$.
  Thus the assignment $X\mapsto \bv(X,D)$ is a contravariant functor
  from the coarse category of coarse spaces to the category of
  $C^*$\nbd{}algebras.
\end{prop}

\begin{proof}
  We identify $\bv(X,D)$ with another $C^*$\nbd{}algebra that is
  evidently functorial for Borel maps.  Let $B_0(X,D\otimes\Comp)$ be
  the $C^*$\nbd{}algebra of bounded Borel functions $X\to
  D\otimes\Comp$ that vanish at infinity.  Let $\Borelvv(X,D)$ consist
  of bounded Borel functions $X\to D\otimes\Comp$ with vanishing
  variation and let $\Borelbv(X,D)\defeq
  \Borelvv(X,D)/B_0(X,D\otimes\Comp)$.  It is evident that
  $B_0(X,D\otimes\Comp)$ and $\Borelvv(X,D)$ and hence $\Borelbv(X,D)$
  are functorial for coarse maps.  Moreover, if $\phi,\phi'\colon X\to
  Y$ are close and $f\in\Borelbv(Y,D)$, then $f\circ\phi-f\circ\phi'$
  vanishes at infinity.  Hence $\phi$ and~$\phi'$ induce the same map
  $\Borelbv(Y,D)\to\Borelbv(X,D)$.

  It is clear that
  $$
  C_0(X,D\otimes\Comp) \subseteq B_0(X,D\otimes\Comp),
  \qquad
  \vv(X,D) \subseteq \Borelvv(X,D).
  $$
  Hence we get an induced $*$\nbd{}homomorphism
  $j\colon \bv(X,D)\to\Borelbv(X,D)$.  We claim that this map is an
  isomorphism.  Once this claim is established, we obtain the desired
  functoriality of $\bv(X,D)$.  Injectivity and surjectivity of~$j$
  are equivalent to
  \begin{align*}
    C_0(X,D\otimes\Comp) &= B_0(X,D\otimes\Comp) \cap \vv(X,D),
    \\
    \Borelbv(X,D) &= \bv(X,D) + B_0(X,D\otimes\Comp),
  \end{align*}
  respectively.  The first equation is evident.  We prove the
  second one.  Let $E\subseteq X\times X$ be an entourage that is
  a neighborhood of the diagonal.  We remarked after
  Definition~\ref{definitionofcoarsestructure} that there exists a
  uniformly bounded open cover $(B_i)$ of~$X$ with $\bigcup B_i\times
  B_i\subseteq E$.  Let~$(\rho_i)$ be a partition of unity subordinate
  to~$(B_i)$ and fix points $x_i\in B_i$.  Take $f\in\Borelbv(X,D)$
  and define
  $$
  Pf(x) \defeq \sum \rho_i (x) f(x_i).
  $$
  It is clear that $Pf$ is continuous.  Since~$f$ has vanishing
  variation, there exists a bounded set $\Sigma\subseteq X$ such that
  $\norm{f(x)-f(y)}<\epsilon$ for $(x,y)\in E$ and $x\notin\Sigma$.
  Hence
  $$
  \norm{(f-Pf)(x)}
  \le \sum \rho_i(x) \norm{f(x)-f(x_i)}
  \le \sum \rho_i(x) \epsilon
  = \epsilon
  $$
  for all $x\notin\Sigma$.  This means that $f-Pf$ vanishes at
  infinity, that is, $f-Pf\in B_0(X,D\otimes\Comp)$.  It follows that
  $Pf\in\vv(X,D)$.  This finishes the proof.
\end{proof}

We are interested in the $\K$\nbd{}theory of the stable Higson corona
$\bv(X,D)$.  We can identify $D\otimes\Comp$ with the subalgebra of
constant functions in $\bv(X,D)$.  It is often convenient to neglect
the part of the $\K$\nbd{}theory that arises from this embedding.
This is the purpose of the following definition.

\begin{defn}  \label{definitionofreducedktheory}
  Let~$X$ be an unbounded coarse space and let~$D$ be a
  $C^*$\nbd{}algebra.  The \emph{reduced} $\K$\nbd{}theory of
  $\vv(X,D)$ and $\bv(X,D)$ is defined by
  \begin{align*}
    \tilde{\K}_*\bigl(\vv(X,D)\bigr)
    & \defeq \K_*\bigl(\vv(X,D)\bigr)/
    \mathrm{range}\bigl[\K_*(D\otimes\Comp)
    \to \K_*\bigl(\vv(X,D)\bigr)\bigr],
    \\
    \tilde{\K}_*\bigl(\bv(X,D)\bigr)
    & \defeq
    \K_*\bigl(\bv(X,D)\bigr)/
    \mathrm{range}\bigl[\K_*(D\otimes\Comp)
    \to \K_*\bigl(\bv(X,D)\bigr)\bigr].
  \end{align*}
\end{defn}

\begin{rmk}  \label{rem:bounded_coarse_space}
  If~$X$ is a bounded coarse space, then
  $\vv(X,D)=C(X,D\otimes\Comp)$ and $\bv(X,D)=0$.  In this case, the
  above definition of $\tilde\K_*\bigl(\bv(X,D)\bigr)$ is not
  appropriate and many things obviously fail.  In order to get true
  results in this trivial case as well, we should define
  $\tilde\K_*(\bv(X,D))\defeq \K_*(\bvr(X,D))$ using the
  $C^*$\nbd{}algebra $\bvr(X,D)$ introduced in
  Definition~\ref{def:vvr} below.
\end{rmk}

\begin{lemma}  \label{injectivityofinclusions}
  Let~$X$ be an unbounded coarse space.  Then the inclusions
  $D\otimes\Comp\to\vv(X,D)$ and $D\otimes\Comp\to\vv(X,D)\to\bv(X,D)$
  induce injective maps in $\K$\nbd{}theory.
\end{lemma}

\begin{proof}
  Let $j\colon C_0(X,D\otimes\Comp)\to\vv(X,D)$ and $\bar{\iota}\colon
  D\otimes\Comp\to\vv(X,D)$ be the inclusions, let $\pi\colon
  \vv(X,D)\to\bv(X,D)$ be the quotient map, and let
  $\iota\defeq\pi\circ\bar\iota\colon D\otimes\Comp\to\bv(X,D)$.  To
  see that $\bar{\iota}_*$ is injective, consider the evaluation map
  $\ev_x\colon \vv(X,D)\to D\otimes\Comp$ for $x\in X$.  This map
  splits~$\bar{\iota}$, from which the assertion follows.
  
  To check that~$\iota_*$ is injective, choose $a\in\K_*(D)$ with $\pi_*
  \bar{\iota}_*(a)=\iota_*(a)=0$.  Hence $\bar{\iota}_*(a)=j_*(b)$ for some
  $b\in\K_*(C_0(X,D\otimes\Comp))$ by the $\K$\nbd{}theory long exact
  sequence.  Let $\ev_x^0\colon C_0(X,D\otimes\Comp)\to D\otimes\Comp$
  denote the restriction of the evaluation map to
  $C_0(X,D\otimes\Comp)$.  Since~$X$ is unbounded and $\K$\nbd{}theory
  is compactly supported, we have $(\ev_x^0)_*(b) = 0$ for all $x\in
  X$ outside some compact set.  Then $a=(\ev_x)_* \bar{\iota}_*(a) =
  (\ev_x)_* j_*(b)= (\ev_x^0)_*(b) = 0$, which concludes the proof.
\end{proof}

\begin{rmk}  \label{rmk:reduced_still_exact}
  For any unbounded coarse space~$X$ and any $C^*$\nbd{}algebra~$D$,
  consider the long exact sequence
  $$
  \xymatrix{
    \K_0(C_0(X)\otimes D) \ar[r] &
    \K_0\bigl(\vv(X,D)\bigr) \ar[r]^{\pi_*} &
    \K_0(\bv (X,D)) \ar[d]_\partial
  \\
    \K_1(\bv(X,D)) \ar[u]^\partial &
    \K_1\bigl(\vv(X,D)\bigr) \ar[l]_{\pi_*} &
    \K_1\bigl(C_0(X)\otimes D) \ar[l]
  }
  $$
  associated to the exact sequence of $C^*$\nbd{}algebras
  $$
  0 \to C_0(X)\otimes D\otimes\Comp \to \vv(X,D) \to \bv(X,D)\to 0.
  $$
  By construction, the map $\iota_*\colon \K_*(D\otimes\Comp) \to
  \K_*\bigl(\bv(X,D)\bigr) $ factors through the map $\pi_*\colon
  \K_*\bigl(\vv(X,D)\bigr)\to\K_*\bigl(\bv(X,D)\bigr)$, whence
  $\partial\circ\iota_* = 0$.  Lemma~\ref{injectivityofinclusions}
  shows that we get a long exact sequence
  $$
  \xymatrix{
    \K_0(C_0(X)\otimes D) \ar[r] &
    \tilde{\K}_0\bigl(\vv(X,D)\bigr) \ar[r]^{\pi_*} &
    \tilde{\K}_0(\bv(X,D)) \ar[d]_{\partial}
  \\
    \tilde{\K}_1(\bv(X,D)) \ar[u]^{\partial} &
    \tilde{\K}_1\bigl(\vv(X,D)\bigr) \ar[l]_{\pi_*} &
    \K_1\bigl(C_0(X)\otimes D).  \ar[l]
  }
  $$
\end{rmk}

Finally, we extend the above definitions to the case of
$\sigma$\nbd{}coarse spaces.  This is necessary to construct the
coarse co\nbd{}assembly map in the next section.

Let $\X = \bigcup X_n$ be a $\sigma$\nbd{}coarse space and let~$D$ be
a $C^*$\nbd{}algebra.  We let
\begin{align*}
  C_0(\X,D) &\defeq
  \{f\colon \X\to D \mid
  \text{$f|_{X_n}\in C_0(X_n,D)$ for all $n\in\N$}\},
  \\
  \vv(\X,D) &\defeq
  \{f\colon \X\to D\otimes\Comp \mid
  \text{$f|_{X_n}\in \vv(X_n,D)$ for all $n\in\N$}\}.
\end{align*}
Both $C_0(\X,D)$ and $\vv(\X,D)$ are $\sigma$\nbd{}$C^*$\brd{}algebras
in the terminology of~\cite{Phillips} with respect to the sequence of
$C^*$\nbd{}seminorms
$$
\norm{f}_n \defeq \sup \{\norm{f(x)} \mid x \in X_n\}.
$$
We evidently have
$$
C_0(\X,D) = \varprojlim C_0(X_n,D\otimes\Comp),
\qquad
\vv(\X,D) = \varprojlim \vv(X_n,D\otimes\Comp),
$$
where $\varprojlim$ denotes the projective limit in the category of
$\sigma$\nbd{}$C^*$\brd{}algebras.

Recall that we assumed the maps $X_n\to X_{n+1}$ to be coarse
equivalences.  Proposition~\ref{functorialityofcoronaalgebra} implies
that the induced maps $\bv(X_{n+1},D)\to \bv(X_n,D)$ are
$*$\nbd{}isomorphisms.  Hence the inverse limit
$$
\bv(\X,D) \defeq \varprojlim \bv(X_n,D)
$$
is again a $C^*$\nbd{}algebra: it is isomorphic to $\bv(X_m,D)$ for
any $m\in\N$.  The following lemma asserts that we also have a natural
isomorphism
$$
\bv(\X,D) \cong \vv(\X,D)/C_0(\X,D\otimes\Comp).
$$

\begin{lemma}  \label{exactness}
  The sequence of $\sigma$\nbd{}$C^*$\brd{}algebras
  $$
  0 \to C_0(\X,D) \to \vv(\X,D) \to \bv(\X,D) \to 0
  $$
  is exact.
\end{lemma}

\begin{proof}[Proof of Lemma~\ref{exactness}]
  The maps $\alpha_n\colon C_0(X_{n+1},D\otimes \Comp) \to
  C_0(X_n,D\otimes\Comp)$ associated to the inclusions $X_n\subseteq
  X_{n+1}$ are clearly surjective.  The maps $\gamma_n\colon
  \bv(X_{n+1},D)\to\bv(X_n,D)$ are surjective because they are
  isomorphisms.  The Snake Lemma of homological algebra provides us
  with a long exact sequence
  $$
  \cdots\to\Coker\alpha_n \to \Coker\beta_n\to\Coker \gamma_n\to0.
  $$
  Hence the maps $\beta_n\colon\vv(X_{n+1},D)\to\vv(X_n,D)$ are
  surjective as well.  Now the assertion follows from the following
  lemma from~\cite{Phillips}.
\end{proof}

\begin{lemma}[\cite{Phillips}]  \label{lem:Phillips}
  Suppose that $\alpha_n\colon A_{n+1} \to A_n$, $n\in\N$, is a
  projective system of $C^*$\nbd{}algebras with surjective
  maps~$\alpha_n$ for all~$n$.  Let~$J_n$ be ideals in~$A_n$ such that
  the restriction of~$\alpha_n$ to $J_{n+1}$ maps $J_{n+1}$
  surjectively onto~$J_n$.  Then
  $$
  0 \to \varprojlim J_n \to \varprojlim A_n \to \varprojlim A_n/J_n
  \to 0
  $$
  is an exact sequence of $\sigma$\nbd{}$C^*$\brd{}algebras.
\end{lemma}

\section{The coarse co-assembly map}
\label{sec:coarse_coassembly}

We first define the coarse $\K$\nbd{}theory of~$X$ with
coefficients~$D$, which is the target of the coarse co\nbd{}assembly
map.  Its definition is based on the Rips complex construction of
Example~\ref{ex:metric_rips}.  We reformulate it in terms of
entourages and check carefully that we obtain a $\sigma$\nbd{}coarse
space.  We require the coarse structure to be countably generated.
Otherwise the construction below gives an uncountable system of coarse
spaces, which we prefer to avoid.  We begin with the case where~$X$ is
discrete.

We fix an increasing sequence $(E_n)$ of entourages such that any
entourage is contained in some~$E_n$.  We assume that $E_0=\Delta_X$
is the diagonal.  Let~$\rips_X$ be the set of probability measures
on~$X$ with finite support.  This is a simplicial complex whose
vertices are the Dirac measures on~$X$.  We give it the corresponding
topology.  Hence locally finite subcomplexes of~$\rips_X$ are locally
compact.  Let
$$
\ripsn_n\defeq \{\mu\in\rips_X \mid
\supp \mu \times \supp \mu \subseteq E_n\}.
$$
In particular, $\ripsn_0\cong X$.  We have $\bigcup \ripsn_n =
\rips_X$ because any finite subset of $X\times X$ is contained
in~$E_n$ for some~$n$.  Each~$\ripsn_n$ is a locally finite subcomplex
of~$\rips_X$ and hence a locally compact topological space because
bounded subsets of~$X$ are finite.  We give~$\rips_X$ and its
subspaces~$\ripsn_n$ the coarse structure~$\mathcal{E}_n$ that is
generated by the sequence of entourages
$$
\{(\mu,\nu) \mid \supp \mu \times \supp \nu \subset E_m\},
\qquad m\in\N.
$$
The embeddings $X\cong \ripsn_0\to\ripsn_n$ are coarse equivalences
for all $n\in\N$.  Thus~$\rips_X$ is a $\sigma$\nbd{}coarse space.

The $\K$\nbd{}theory of the $\sigma$\nbd{}$C^*$\nbd{}algebra
$C_0(\rips_X,D)$ is going to be the coarse $\K$\nbd{}theory of~$X$.
In order to extend this definition to non-discrete coarse spaces, we
must show that it is functorial on the coarse category of coarse
spaces.

We first observe that the $\sigma$\nbd{}$C^*$\brd{}algebra
$C_0(\rips_X,D)$ does not depend on the choice of the generating
sequence~$(E_n)$.  A function $f\colon \rips_X\to D$ belongs to
$C_0(\rips_X,D)$ if and only if its restriction to $\ripsn_n$ is~$C_0$
for all $n\in\N$.  If $E'\subseteq X\times X$ is any entourage, then
$E'\subseteq E_n$ for some $n\in\N$.  If we define
$\ripsn_{E'}(X)\subseteq \rips_X$ in the evident fashion, we obtain a
subcomplex of $\ripsn_n(X)$.  Thus the restriction of~$f$ to
$\ripsn_{E'}(X)$ is~$C_0$ for all entourages~$E'$.  Conversely, this
condition implies easily that $f\in C_0(\rips_X,D)$.  Hence we can
describe $C_0(\rips_X,D)$ without using the generating sequence~$E_n$.
Similar arguments apply to $\vv(\rips_X,D)$ and, of course, to
$\bv(\rips_X,D)$.  Actually, up to an appropriate notion of
isomorphism of inductive systems, the $\sigma$\nbd{}coarse
space~$\rips_X$ is independent of the choice of~$(E_n)$.

To discuss the functoriality of~$\rips_X$, we define morphisms between
$\sigma$\nbd{}coarse spaces.  Let $\X=\bigcup X_n$ and $\Y=\bigcup
Y_n$ be $\sigma$\nbd{}coarse spaces.  Let $f\colon \X\to\Y$ be a map
with the property that for any $m\in\N$ there is $n=n(m)\in\N$ with
$f(X_m)\subseteq Y_n$.  We say that~$f$ is Borel, continuous, or
coarse, respectively, if the restrictions $f|_{X_m}\colon X_m\to
Y_{n(m)}$ have this property for all $m\in\N$.  Here the choice of
$n(m)$ is irrelevant because~$Y_n$ is a subspace of~$Y_{n'}$ for all
$n\le n'$.  Similarly, two coarse maps $\X\to\Y$ are called close if
their restrictions to~$X_m$ are close for all $m\in\N$.  It is clear
that $C_0(\X,D)$ and $\vv(\X,D)$ are functorial for continuous coarse
maps.

\begin{lemma}  \label{lem:rips_functorial}
  Let~$X$ and $Y$ be discrete, countably generated coarse spaces.
  Then a coarse map $X \to Y$ induces a continuous coarse map
  $\rips_X\to\rips_Y$.

  Let $\X$ be a $\sigma$\nbd{}coarse space and let
  $\phi_0,\phi_1\colon\X\to\rips_Y$ be two continuous coarse maps that
  are close.  Then there exists a homotopy $\Phi\colon
  \X\times[0,1]\to\rips_Y$ between $\phi_0$ and~$\phi_1$ that is close
  to the constant homotopy $(x,t)\mapsto \phi_0(x)$.  The
  homotopy~$\Phi$ induces a homotopy $C_0(\rips_Y)\to C([0,1])\otimes
  C_0(\X)$.
\end{lemma}

\begin{proof}
  A coarse map $\phi\colon X\to Y$ induces a map $\phi_*\colon
  \rips_X\to\rips_Y$ by pushing forward probability measures.  It is
  easy to see that~$\phi_*$ is continuous and coarse.
  
  We want to define $\Phi(x,t)\defeq (1-t)\phi_0(x)+t\phi_1(x)$.  It
  is clear that $\Phi(x,t)$ is a probability measure on~$Y$ with
  finite support for all $(x,t)\in\X\times [0,1]$.  We claim
  that~$\Phi$ has the required properties.  Fix $n\in\N$ and an
  entourage $E\subseteq X_n\times X_n$.  Since $\phi$ and~$\phi'$ are
  close, there is an entourage $E'\subseteq Y\times Y$ such that
  $$
  \phi_0\times\phi_1(E) \subseteq
  \{(\mu,\nu) \mid \supp \mu \times \supp \nu \subseteq E'\}.
  $$
  Let $E''\defeq E'\cup (E'\circ (E')^t)$.  Since $\supp \Phi(x,t)
  \subseteq \supp \phi_0(x)\cup\supp \phi_1(x)$, we obtain
  $\Phi(x,t)\in\ripsn_{E''}$ and $\supp \phi_0(x) \times \supp
  \Phi(x,t) \subseteq E''$ for all $x\in X_n$, $t\in[0,1]$.  That is,
  $\Phi$ is a coarse map that is close to the map $(x,t)\mapsto
  \phi_0(x)$.  Continuity is easy to check.  We also get an induced
  homotopy for the associated $\sigma$\brd{}$C^*$\nbd{}algebras
  because $C_0(\X,D)\otimes C([0,1]) \cong C_0(\X\times[0,1],D)$.
\end{proof}

Since $\K$\nbd{}theory for $\sigma$\brd{}$C^*$\nbd{}algebras is still
homotopy invariant, we obtain:

\begin{cor}  \label{cor:coarse_K_functorial}
  The assignment $X\mapsto \K_*\bigl(C_0(\rips_X,D)\bigr)$ is a
  functor from the coarse category of discrete, countably generated
  coarse spaces to the category of $\Ztwo$\nbd{}graded Abelian groups.
\end{cor}

\begin{defn}
  Let~$X$ be a countably generated coarse space and let~$D$ be a
  $C^*$\nbd{}algebra.  Let $Z\subseteq X$ be a countably generated,
  discrete coarse space that is coarsely equivalent to~$X$.  This
  exists by Lemma~\ref{lem:discretization}.  We let
  $$
  \KX^*(X,D)\defeq
  \K_*\bigl(C_0(\rips_Z,D)\bigr)
  $$
  and call this the \emph{coarse $\K$\nbd{}theory of~$X$ with
  coefficients~$D$}.
\end{defn}

\begin{note}
  When $D=\C$ is trivial we simply write $\KX^*(X)\defeq\KX^*(X,\C)$
  and refer to this as the \emph{coarse $\K$\nbd{}theory of~$X$}.
\end{note}

By construction, the discrete coarse space~$Z$ is uniquely determined
up to coarse equivalence.  Hence
Corollary~\ref{cor:coarse_K_functorial} yields that $\KX^*(X,D)$ is
independent of the choice of~$Z$ and is a functor from the coarse
category of countably generated coarse spaces to the category of
$\Ztwo$\nbd{}graded Abelian groups.  Since the homotopy type
of~$\rips_Z$ is independent of the choice of~$Z$, we also
write~$\rips_X$ for this space.

\begin{rmk}  \label{rmk:Milnor_sequence}
  Phillips shows in~\cite{Phillips} that $\K$\nbd{}theory for
  $\sigma$\nbd{}$C^*$\brd{}algebras can be computed by a Milnor
  $\varprojlim^1$\brd{}sequence.  In our case, we obtain a short exact
  sequence
  $$
  0 \to
  \varprojlim\nolimits^1 \K_{*+1}(C_0(\ripsn_n(Z),D)) \to
  \KX^*(X,D) \to
  \varprojlim \K_*(C_0(\ripsn_n(Z),D)) \to 0.
  $$
  For $D=\C$, this becomes a short exact sequence
  $$
  0 \to
  \varprojlim\nolimits^1 \K^{*+1}(\ripsn_n(Z)) \to
  \KX^*(X) \to
  \varprojlim \K^*(\ripsn_n(Z)) \to 0.
  $$
\end{rmk}

We are now in a position to define our coarse co\brd{}assembly map.
Let~$D$ be a $C^*$\nbd{}algebra and let~$X$ be a coarse space.  By
Lemma~\ref{exactness} the sequence
\begin{equation}  \label{fundamentalexactsequence}
  0\to C_0(\rips_X,D\otimes\Comp) \to \vv(\rips_X,D) \to
  \bv(\rips_X,D) \cong \bv(X,D) \to 0
\end{equation}
is exact.  In~\cite{Phillips} it is shown that an exact sequence of
$\sigma$\nbd{}$C^*$\brd{}algebras induces a long exact sequence in
$\K$\nbd{}theory.  As in Remark~\ref{rmk:reduced_still_exact}, one
shows that this remains exact if we use reduced $\K$\nbd{}theory
everywhere (and assume~$X$ to be unbounded).

\begin{defn}  \label{definitionofthecoarsecoassemblymap}
  Let~$X$ be a countably generated, unbounded coarse space and let~$D$
  be a $C^*$\nbd{}algebra.  The \emph{coarse co\brd{}assembly map
    for~$X$ with coefficients~$D$} is the map
  $$
  \mu^*_{X,D}\colon \tilde\K_{*+1}\bigl(\bv(X,D)\bigr)\to\KX^*(X,D)
  $$
  that is obtained from the connecting map of the exact
  sequence~\eqref{fundamentalexactsequence}.
\end{defn}

We conclude this section by noting that the coarse $\K$\nbd{}theory
of~$X$ is equal to the usual $\K$\nbd{}theory of~$X$ if~$X$ is a
uniformly contractible metric space of bounded geometry.  An analogous
assertion holds for the coarse $\K$\nbd{}homology.

\begin{defn}  \label{definitionofuniformlycontractible}
  A metric space $(X,d)$ is \emph{uniformly contractible} if for every
  $R>0$ there exists $S\ge R$ such that for any $x\in X$, the
  inclusion $B_R(x) \to B_S(x)$ is nullhomotopic.
\end{defn}

\begin{thm}  \label{identificationofcoarsektheory}
  Let~$X$ be a uniformly contractible metric space of bounded
  geometry.  There exists a canonical isomorphism $\KX^*(X,D)\cong
  \K_*\bigl(C_0(X,D)\bigr)$ making the following diagram commute:
  \begin{displaymath}
    \xymatrix@C=-1em{
      {\K_{*+1}\bigl(\bv(X,D)\bigr)} \ar[rr]^-{\mu_{X,D}^*}
      \ar[dr]_{\partial} & &
      {\KX^*(X,D)} \\
      & {\K_*\bigl( C_0 (X,D)\bigr).}  \ar[ur]_{\cong}
      &
    }
  \end{displaymath}
  Here~$\partial$ is the boundary map associated to the exact sequence
  of $C^*$\nbd{}algebras
  $$
  0 \to C_0(X,D\otimes \Comp) \to \vv(X,D) \to \bv(X,D) \to 0.
  $$
\end{thm}

The proof uses the following lemma. 

\begin{lemma}  \label{homotopylemma}
  Let~$X$ be a uniformly contractible metric space of bounded
  geometry, and let $\phi\colon X\to X$ be a continuous coarse map
  which is close to the identity map $X\to X$.  Then $\phi$ and
  $\mathrm{id}$ are homotopic, and the homotopy $F\colon
  X\times[0,1]\to X$ can be chosen to be close to the coordinate
  projection $(x,t)\mapsto x$.
\end{lemma}

\begin{proof}
  Choose a uniformly bounded open cover $(U_i)_{i\in I}$
  of~$X$.  For $\Sigma\subseteq I$ let $U_\Sigma = \bigcap_{i\in\Sigma}
  U_i$ and let $\Delta_\Sigma$ be the simplex $\Delta_\Sigma \defeq
  \{(x_i)\in[0,1]^\Sigma \mid \sum x_i=1\}$.  By the bounded geometry
  assumption we may choose the cover so that $U_\Sigma=\emptyset$
  whenever $\abs\Sigma > N$ for some $N\in\N$.  By induction on
  $n=\abs\Sigma\ge1$, we construct continuous maps $H_\Sigma\colon
  U_\Sigma\times \Delta_\Sigma\times [0,1]\to X$ such that
  \begin{enumerate}[(1)]
  \item $H_\Sigma(x,p,0)=x$ and $H_\Sigma(x,p,1)=\phi(x)$ for all
    $x\in U_\Sigma$, $p\in\Delta_\Sigma$;

  \item $H_\Sigma(x,p,t)=H_{\supp(p)}(x,p,t)$ for all $x\in U_\Sigma$,
    $p\in\partial\Delta_\Sigma$, $t\in[0,1]$;

  \item $d(H_\Sigma(x,p,t),x) \le C_n$ for all $x,p,t$ for some
    constant $C_n\ge 0$.

  \end{enumerate}
  
  If $\abs\Sigma=1$, then $\Sigma=\{i\}$ for some $i\in I$, $U_\Sigma
  = U_i$ and $\Delta_\Sigma$ is a point, denote it~$\star$.  The map
  $H_\Sigma$ must satisfy $H_\Sigma(x,\star,0) = x$ and
  $H_\Sigma(x,\star,1)=\phi(x)$ for all $x\in U_i$.  By the uniform
  contractibility, we can extend this to a continuous map on
  $U_i\times[0,1]$ with the required properties.  Now assume
  that~$H_\Sigma$ has been defined for $\abs\Sigma<n$ and take
  $\Sigma\subseteq I$ with $\abs\Sigma=n$.  The previous induction
  step and our requirements determine~$H_\Sigma$ on $(U_\Sigma \times
  \partial \Delta_\Sigma \times [0,1]) \cup (U_\Sigma \times
  \Delta_\Sigma \times \partial [0,1])$.  The uniform contractibility
  assumption of~$X$ allows us to extend this to $U_\Sigma \times
  \Delta_\Sigma \times [0,1]$ as required.
  
  Finally, we choose a partition of unity $(\rho_i)_{i\in I}$
  subordinate to the cover $(U_i)_{i \in I}$ and define $F\colon
  X\times [0,1] \to X$ as follows.  For $x\in X$, let $\Sigma\defeq
  \{i\in I \mid \rho_i(x) \neq 0\}$ and $F(x,t)\defeq
  H_\Sigma(x,(\rho_i(x))_{i\in\Sigma},t)$.  This defines a continuous
  map that is close to the identity map because $U_\Sigma=\emptyset$
  for $\abs\Sigma>N$.
\end{proof}

\begin{proof}[Proof of Theorem~\ref{identificationofcoarsektheory}]
  Let $Z\subseteq X$ be a discrete subspace coarsely equivalent to~$X$
  as in Lemma~\ref{lem:discretization}.  For sufficiently large~$r$,
  the balls of radius~$r$ centered at the points of~$Z$ cover~$X$.
  Let $\ripsn_r(Z)$ be the Rips complex with parameter~$r$ as in
  Example~\ref{ex:metric_rips}.  The natural maps $Z\to X$ and $Z\to
  \ripsn_r(Z)$ are coarse equivalences.  Hence we obtain canonical
  coarse equivalences $X\to \ripsn_r(Z)$ and $\ripsn_r(Z)\to X$.  We
  want to represent these morphisms in the coarse category by
  continuous coarse maps $F\colon X\to \ripsn_r(Z)$ and $G\colon
  \ripsn_r(Z)\to X$.
  
  Let $(\rho_i)$ be a partition of unity subordinate to the cover
  of~$X$ by $r$\nbd{}balls centered at the points of~$Z$.  Choose
  $x_i\in Z$ close to $\supp\rho_i$ and define
  $$
  F\colon X\to \ripsn_r(Z)\subseteq \rips_Z,
  \qquad
  F(x)\defeq \sum_i \rho_i(x) \delta_{x_i}.
  $$
  This is a continuous coarse map, and its restriction to~$Z$ is close
  to the standard map $Z\to \ripsn_r(Z)$ as desired.  Notice that this
  map exists for any countably generated coarse space~$X$.
  
  We define the maps $G\colon \ripsn_r(Z)\to X$ for any $r\ge0$ by
  induction on skeleta.  On the $0$\nbd{}skeleton~$Z$, we let~$G$ be
  the inclusion $Z\to X$.  Suppose that~$G$ has already been defined
  on the $n-1$\brd{}skeleton and let~$\sigma$ be an $n$\nbd{}cell.
  Then the vertices of~$\sigma$ constitute a subset of $Z\subseteq X$
  of diameter at most~$r$.  By our induction assumption, $G$ maps the
  boundary of~$\sigma$ to a subset of~$X$ of diameter at most
  $C_{n-1}(r)$ for some constant depending only on $r$ and $n-1$.  By
  uniform contractibility, we can extend~$G$ to a map $\sigma\to X$ in
  such a way that $G(\sigma)$ has diameter at most $C_n(r)$ for some
  constant $C_n(r)$.  Proceeding in this fashion, we construct a
  continuous coarse map $G\colon \ripsn_r(Z)\to X$ whose restriction
  to~$Z$ is the inclusion map.

  The compositions $F\circ G$ and $G\circ F$ are continuous coarse
  maps which are close to the identity maps on $\ripsn_r(Z)$ and $X$
  respectively.  Lemmas \ref{lem:rips_functorial}
  and~\ref{homotopylemma} yield that $F\circ G$ and $G\circ F$ are
  homotopic to the identity maps on $\rips_Z$ and~$X$, respectively.
  Therefore, we get an isomorphism $\K_*(C_0(X,D))\cong
  \K_*(C_0(\rips_Z,D))$ as desired.  Since $F$ and~$G$ extend to
  $*$\nbd{}homomorphisms between $\vv(X,D)$ and $\vv(\rips_X,D)$, the
  naturality of the boundary map in $\K$\nbd{}theory yields the
  commutative diagram in the statement of the theorem.
\end{proof}

\section{A first vanishing theorem}
\label{sec:vanishing}

We now calculate an example that illustrates the distinction between
the stable and unstable Higson coronas.  We begin by recalling the
following result of~\cite{Keesling}.

\begin{prop}
  Let $X=[0,\infty)$ be the ray with its (Euclidean) metric coarse
  structure.  Then the reduced $\K$\nbd{}theory of the Higson
  compactification~$\eta X$ of~$X$ is uncountable.
\end{prop}

This implies by the Five Lemma that the same is true for the Higson
corona $\partial_\eta X$ of~$X$.  In contrast, we show that the
reduced $\K$\nbd{}theory of the stable Higson corona of~$X$ is
trivial.  Since it involves no additional effort, we show the
following more general result.

\begin{thm}  \label{firstvanishingtheorem}
  Let~$Y$ be an arbitrary coarse space and~$D$ a $C^*$\nbd{}algebra.
  Let the ray $[0,\infty)$ be given its Euclidean coarse structure
  and let $X=Y\times[0,\infty)$ with the product coarse structure.  Then
  $\tilde{\K}_*\bigl(\bv(X,D)\bigr)=0$ for $*=0,1$.
\end{thm}

\begin{rmk}
  This result is consistent with the analogous assertion
  $\K_*\bigl(C^*(X)\bigr)=0$ for the Roe $C^*$\nbd{}algebras of such
  spaces.
\end{rmk}

For the purposes of this computation and for many others, it turns out
to be much easier to work not with the reduced $\K$\nbd{}theory of the
algebras $\vv$ and~$\bv$, but rather with the ordinary
$\K$\nbd{}theory of modified (or reduced) versions of these algebras.
Thus we introduce the following definition.  If~$D$ is a
$C^*$\nbd{}algebra, we let $\Mult(D)$ be its multiplier algebra and
$\Mult^s(D)$ be the multiplier algebra of $D\otimes\Comp$.  We also
define $\Calkin^s(D)\defeq \Mult^s(D)/D\otimes\Comp$.

\begin{defn}  \label{def:vvr}
  Let~$X$ be a coarse space and let~$D$ be a $C^*$\nbd{}algebra.  We
  let $\vvr(X,D)$ be the $C^*$\nbd{}algebra of bounded continuous
  functions of vanishing variation $f\colon X\to\Mult^s(D)$ such
  that $f(x)-f(y)\in D\otimes\Comp$ for all $x,y\in X$.  We let
  $\bvr(X,D)$ be the $C^*$\nbd{}algebra
  $\vvr(X,D)/C_0(X,D\otimes\Comp)$.
\end{defn}

\begin{prop}  \label{identificationofreducedktheory}
  For every unbounded metric space~$X$ and every
  $C^*$\nbd{}algebra~$D$, we have natural isomorphisms
  \begin{align*}
    \K_*\bigl(\vvr(X,D)\bigr) &\cong \tilde{\K}_*\bigl(\vv(X,D)\bigr),
    \\
    \K_*\bigl(\bvr(X,D)\bigr) &\cong \tilde{\K}_*\bigl(\bv(X,D)\bigr).
  \end{align*}
\end{prop}

\begin{proof}
  Choose any point $x\in X$, and consider the composition
  $$
  \vvr(X,D)\to \Mult^s(D)\to \Calkin^s(D),
  $$
  where the first map is evaluation at $x\in X$ and the second is
  the quotient map.  This map is surjective and its kernel is the
  unreduced algebra $\vv(X,D)$.  Therefore, it descends to a map on
  the quotient $\bvr(X,D)$.  Applying the Bott Periodicity
  isomorphisms $\K_*\bigl(\Calkin^s(D)\bigr) \cong \K_{*+1} (D)$,
  we obtain a long exact sequence
  $$
  \xymatrix{
    \K_0\bigl(\vv(X,D)\bigr) \ar[r] &
    \K_0\bigl(\vvr(X,D)\bigr) \ar[r] &
    \K_1(D) \ar[d]^{\bar\iota_*}
    \\
    \K_0(D) \ar[u]^{\bar\iota_*} &
    \K_1\bigl(\vvr(X,D)\bigr) \ar[l] &
    \K_1\bigl(\vv(X,D)). \ar[l]
  }
  $$
  One can show that the boundary maps are induced by the inclusion
  $\bar\iota\colon D\otimes\Comp\to\vv(X,D)$.  Since the vertical maps
  are injective by Lemma~\ref{injectivityofinclusions}, we obtain two
  exact sequences
  $$
  0\to \K_*(D)\overset{\bar\iota_*}\to \K_*\bigl(\vv(X,D)\bigr)\to
  \K_*\bigl(\vvr(X,D)\bigr) \to 0
  $$
  for $*=0,1$.  This proves the first assertion.  The second one is
  proved in the same fashion.
\end{proof}

\begin{proof}[Proof of Theorem~\ref{firstvanishingtheorem}]
  We may replace $Y\times[0,\infty)$ by the coarsely equivalent space
  $Y\times\N$.  Thus we let $X\defeq Y\times\N$ with its product
  coarse structure.  By Lemma~\ref{identificationofreducedktheory} it
  suffices to calculate the $\K$\nbd{}theory of the algebras $\vvr
  (X,D)$.  Implicit in the definition of $\vvr(X,D)$ is a Hilbert
  space.  Let us integrate this Hilbert space temporarily into our
  notation by denoting $\vvr(X,D)$, built on the Hilbert space~$V$, by
  $\vvr(X,D,V)$.
  
  Now fix a Hilbert space~$H$, and let $\tilde{H}\defeq H\oplus
  H\oplus\cdots$.  The $C^*$\nbd{}algebras $\vv(X,D,H)$ and
  $\vv(X,D,\tilde{H})$ are (non-canonically) isomorphic.  The
  inclusion $H\to\tilde{H}$ as the $n$th summand induces a
  $*$\nbd{}homomorphism $i_n\colon \vvr(X,D,H) \to
  \vvr(X,D,\tilde{H})$ of the form
  $$
  i_n(f)(x) =
  0\oplus\dots\oplus 0\oplus f(x) \oplus0\oplus 0\oplus \cdots \in
  \Bound (D\otimes\tilde{H})
  $$
  for $n\in\N$.  Standard arguments yield that the maps $(i_n)_*$
  all induce isomorphisms on $\K$\nbd{}theory and that $(i_n)_* =
  (i_m)_*$ for any $n,m\in\N$.  Define a $*$\nbd{}homomorphism
  \begin{displaymath}
    S\colon \vvr(X,D,H)\to \vvr(X,D,H),
    \qquad Sf(y,n) \defeq f(y,n+1).
  \end{displaymath}
  The variation condition implies that $\norm{Sf(y,n)-f(y,n)} =
  \norm{f(y,n+1)-f(y,n)}\to 0$ for $(y,n)\to\infty$ for all
  $f\in\vvr(X,D,H)$.  That is, $Sf-f\in C_0(X,D\otimes\Comp)$.
  Hence~$S$ induces the identity map $\bvr(X,D,H)\to\bvr(X,D,H)$.
  
  We claim that $\tilde{S}f \defeq \oplus_{n=0}^{\infty} i_n(S^nf)$,
  that is,
  $$
  \tilde{S}f(y,n) = f(y,n) \oplus f(y,n+1) \oplus f(y,n+2) \oplus
  \ldots,
  $$
  defines a $*$\nbd{}homomorphism $\tilde{S}\colon
  \vvr(X,D,H)\to\vvr(X,D,\tilde{H})$.  Since~$f$ is bounded,
  $\tilde{S}f(x)$ is a bounded operator for all $x\in X$.  We claim
  that $\tilde{S}f(x) - \tilde{S}f(x')$ is compact for all $x=(y,n)$,
  $x'=(y',n')$ in $Y\times\N$.  The $k$th direct summand of
  $\tilde{S}f(x) - \tilde{S}f(x')$ is given by $f(y,n+k)-f(y',n'+k)$
  and hence lies in $D\otimes\Comp(H)$.  The sequence
  $\bigl((y,n+k),(y',n'+k)\bigr)_{k\in\N}$ in $X\times X$ lies in an
  entourage and converges to~$\infty$.  Hence the vanishing variation
  of~$f$ implies
  $$
  \lim_{k\to\infty} \norm{f(y,n+k)-f(y',n'+k)} = 0.
  $$
  Consequently, $\tilde{S}f(x)-\tilde{S}f(x')$ is a compact
  operator as claimed.  The same reasoning shows that~$\tilde{S}f$
  satisfies the variation condition.  Thus
  $\tilde{S}f\in\vvr(X,D,\tilde{H})$.
  
  Now let $f\in\bvr(X,D,H)$ represent a class~$[f]$ in either
  $\K_0\bigl(\bvr(X,D,H)\bigr)$ or $\K_1\bigl(\bvr(X,D,H)\bigr)$.  We
  have to show that $[f]=0$.  Recall that~$S$ represents the identity
  map on $\bvr(X,D,H)$ and that $(i_n)_* = (i_m)_*$ for all $n,m$.
  Hence
  $$
  [\tilde{S}f] =
  \Bigl[ \oplus_{n=0}^\infty i_n\circ S^n(f)\Bigr] =
  (i_0)_*[f] + [\tilde{S}f].
  $$
  Since~$(i_0)_*$ is an isomorphism, we get $[f]=0$ as desired.
\end{proof}

\begin{rmk}
  The proof of Theorem~\ref{firstvanishingtheorem} exhibits the
  difference between $\bv(X)$ and $C(\eta X,\Comp)$.  If $f\colon \eta
  X\to\Comp$ is a continuous function, then~$f$ must satisfy the
  variation condition and, in addition, have compact range in~$\Comp$.
  The function $\tilde{S}f$ need not have compact range and therefore
  can only be formed in the larger algebra $\bv(X)$.
\end{rmk}

\begin{rmk}  \label{mayervietorisargument}
  Theorem~\ref{firstvanishingtheorem}, a Mayer-Vietoris argument and
  induction can be used to prove that the reduced $\K$\nbd{}theory for
  $\bv(\R^n)$ is given by
  $$
  \tilde{\K}_i\bigl(\bv(\R^n)\bigr) \cong
  \begin{cases}
    \Z & \text{if $i=n-1$,} \\
    0  & \text{otherwise}.
  \end{cases}
  $$
  We omit the argument because this also follows from our results on
  scalable spaces.  This calculation shows clearly that the algebra
  $\bv(X)$ plays the role, at least $\K$\nbd{}theoretically, of a
  boundary of~$X$.
\end{rmk}

\section{Relationship with the coarse Baum-Connes assembly map}
\label{sec:duality}

Let~$X$ be a coarse space and define $\rips_X = \bigcup
\ripsn_n$ as above.  One can extend $\K$\nbd{}homology and even
bivariant $\KK$\brd{}theory from the category of $C^*$\nbd{}algebras
to the category of $\sigma$\nbd{}$C^*$\brd{}algebras, see
\cites{Weidner, Bonkat}.  For $\K$\brd{}homology, one obtains
\begin{displaymath}
  \K^* (C_0(\rips_X)\bigr) =
  \K^*\bigl(\varprojlim C_0(\ripsn_n)\bigr) \cong
  \varinjlim \K^*\bigl(C_0(\ripsn_n)\bigr) \cong
  \varinjlim \K_*(\ripsn_n).
\end{displaymath}
The latter is, by definition, the \emph{coarse $\K$\nbd{}homology
  of~$X$} (see~\cite{Yu}).  Thus we get a natural isomorphism
$$
\KX_*(X) \cong \K^* (C_0(\rips_X)\bigr).
$$
The canonical pairing between the $\K$\nbd{}theory and
$\K$\nbd{}homology for $\sigma$\nbd{}$C^*$\brd{}algebras specializes
to a natural pairing
$$
\KX^*(X) \times \KX_*(X) \to\Z.
$$

Let $C^*(X)$ be the $C^*$\nbd{}algebra of the coarse space~$X$ (see
\cites{HigsonRoe, Yu}).  The coarse Baum-Connes assembly map for~$X$ is
a map
$$
\mu\colon \KX_*(X) \to \K_*\bigl(C^*(X)\bigr).
$$
The next theorem asserts that this map is dual to our coarse
co\brd{}assembly map
$$
\mu^*\colon \tilde\K_{*+1}\bigl(\bv(X)\bigr) \to \KX^*(X).
$$

\begin{thm}  \label{the:duality}
  Let~$X$ be an unbounded, countably generated coarse space~$X$.  Then
  there exists a natural pairing
  $$
  \tilde\K_{*+1}\bigl(\bv(X)\bigr) \times \K_*\bigl(C^*(X)\bigr)
  \to \Z
  $$
  compatible with the pairing $\KX^*(X)\times\KX_*(X)\to\Z$ in the
  sense that
  $$
  \braket{\mu(x)}{y} = \braket{x}{\mu^*(y)} \qquad
  \text{for all $x\in\KX_*(X)$,
  $y\in\tilde\K_{*+1}\bigl(\bv(X)\bigr)$.}
  $$
\end{thm}

\begin{proof}
  Without loss of generality, we may assume~$X$ discrete.  Let~$W$ be
  a separable Hilbert space, and form the ample Hilbert space $H_X =
  \ell^2(X)\otimes W$ over~$X$ (see~\cite{HigsonRoe}).  We use~$H_X$ to
  construct $C^*(X)$.  Thus $C^*(X)$ becomes the $C^*$\nbd{}subalgebra
  of $\Bound(H_X)$ generated by the $*$\nbd{}algebra of locally
  compact finite propagation operators on~$H_X$.  Let~$V$ be another
  separable Hilbert space and let $\Comp\cong\Comp(V)$ be represented
  on~$V$ in the obvious way.  Let $(e_x)_{x\in X}$ be the canonical
  basis of $\ell^2(X)$.  We represent $\vvr(X)$ on $H_X\otimes V$ by
  the map $f\mapsto M_f$ with
  $$
  M_f (e_x\otimes w\otimes v) = e_x \otimes w \otimes f(x)v,
  $$
  for all $f\in\vvr(X)$, viewed as a function $X\to\Comp(V)$.
  Represent $C^*(X)$ on $H_X\otimes V$ by $T\mapsto T\otimes 1_V$.
  The variation condition on $f\in\vvr(X)$ and the definition of
  finite propagation imply easily that the commutator $[M_f,T\otimes
  1_V]$ is compact for all $f\in\vvr(X)$ and $T\in C^*(X)$.  Hence we
  have defined a $*$\nbd{}homomorphism from $\vvr(X)$ into
  $\curD\bigl(C^*(X)\bigr)$ in the notation
  of~\cite{HigsonRoe}.  If $f\in C_0(X,\Comp)$, then both $M_f\cdot
  (T\otimes 1_V)$ and $(T\otimes 1_V)\cdot M_f$ are compact.  That is,
  $C_0(X,\Comp)$ is mapped to $\curD\bigl(C^*(X) // C^*(X) \bigr)$.
  Hence $\bv(X)$ is mapped to the relative dual
  $$
  \relD\bigl(C^*(X)\bigr) \defeq
  \curD\bigl(C^*(X)\bigr) /\curD\bigl(C^*(X) // C^*(X) \bigr)
  $$
  and we obtain a map
  $\K_*\bigl(\bvr(X)\bigr) \to \K_*\bigl(\relD(C^*(X))\bigr)$.
  For every $C^*$\nbd{}algebra~$A$ regardless of separability there
  is a canonical index pairing
  $$
  \K_{*+1}(A) \times \K_*\bigl(\relD(A)\bigr) \to \Z.
  $$
  Since $\tilde\K(\bv(X))\cong \K(\bvr(X))$, we obtain the required
  pairing.

  We omit the details of the proof that $\mu$ and~$\mu^*$ are
  compatible.  First, one shows that it suffices to look at a fixed
  parameter in the Rips complex construction.  The result then follows
  from the definition of~$\mu$ given in terms of dual algebras
  (see~\cite{HigsonRoe}), our definition of~$\mu^*$, and the axioms
  for a Kasparov product.
\end{proof}

\begin{cor}
  Let~$X$ be a uniformly contractible metric space of bounded
  geometry, endowed with the metric coarse structure.  If the coarse
  co\brd{}assembly map for~$X$ is surjective, then the coarse assembly
  map is rationally injective.
\end{cor}

\begin{proof}
  In this case, we can use~$X$ itself instead of the Rips complex by
  Theorem~\ref{identificationofcoarsektheory}.  Hence the pairing
  between $\KX_*(X)\otimes\Q$ and $\KX^*(X)\otimes\Q$ is
  non\brd{}degenerate.
\end{proof}

In particular, let~$X$ be the universal cover of a compact aspherical
spin manifold.  The surjectivity of~$\mu^*$ for~$X$ implies rational
injectivity of~$\mu$ for~$X$.  This in turn implies that the manifold
does not admit a metric of positive scalar curvature.

A natural question is whether or not rational surjectivity of~$\mu_X$
can be detected by injectivity, or even bijectivity, of~$\mu^*_X$.
There seems, however, little hope for this, as the following example
illustrates.  Let $X_1,X_2,\ldots$ be a sequence of finite metric
spaces and let $X=\bigsqcup X_n$ be their \emph{coarse (uniform)
disjoint union}, whose coarse structure is generated by entourages of
the form $\bigsqcup_{i=1}^\infty E_{R,i}$, where~$E_{R,i}$ is the
entourage of diameter~$R$ in~$X_i$.

\begin{prop}
  Let $(X_n)$ be a sequence of finite metric spaces and let $X =
  \bigsqcup X_n$ be the coarse disjoint union as above.  Then the
  pairing
  $$
  \K_*\bigl(C^*(X)\bigr) \times \tilde\K_{*+1}\bigl(\bv(X)\bigr)
  \to \Z
  $$
  is the zero map.
\end{prop}

That is, it is impossible in this example to detect elements of
$\K_*\bigl(C^*(X)\bigr)$ by pairing them with the $\K$\nbd{}theory
of the stable Higson corona.

\begin{proof}
  We use the ample Hilbert space $H_X=\ell^2(X)\otimes W \cong
  \bigoplus l^2(X_n) \otimes W$ to realize $C^*(X)$.  Let~$T$ be a
  finite propagation operator on~$X$.  Then~$T$ can be represented by
  a block diagonal operator $T = S \oplus T_N \oplus T_{N+1} \oplus
  \cdots$ with operators~$T_i$ on~$X_i$ of uniform finite propagation
  and an operator~$S$ that is supported on the bounded set
  $\bigcup_{i=1}^N X_i$.  Note that~$S$ and the operators~$T_i$ are
  compact.  Since the sum of $\Comp(H_X)$ and $\prod \Comp(H_{X_n})$
  is a $C^*$\nbd{}algebra, it contains $C^*(X)$.  Thus any
  element~$\alpha$ of $\K_0(C^*(X))$ is represented by a block
  diagonal projection~$P$ on~$H_X$ with compact blocks.  To construct
  the pairing between $\K_0\bigl(C^*(X)\bigr)$ and
  $\tilde{\K}_1\bigl(\bv(X)\bigr)$, one forms the Hilbert space
  $H_X\otimes V\cong \bigoplus l^2(X_n) \otimes W\otimes V$ as before.

  Let $f\in\bvr(X)$ be a unitary element representing an element of
  $\tilde\K_1(\bv(X))$ and let~$\bar{f}$ be a lifting of~$f$ to an
  element of $\vvr(X)$.  Hence $\bar{f}\bar{f}^*-1$ and
  $\bar{f}^*\bar{f}-1$ lie in $C_0(X,\Comp)$.  The pairing
  $\braket{\alpha}{[f]}$ is given by the index of the Fredholm
  operator $PM_{\bar{f}}P + 1-P$.  This operator is also block
  diagonal, and its blocks are compact perturbations of~$1$.  Hence
  the index vanishes.  The same argument works for the pairing between
  $\K_1(C^*(X))$ and $\tilde\K_0(\bv(X))$.
\end{proof}

\begin{rmk}
  Let~$E$ be a coarse disjoint union of graphs~$E_n$.  Let~$\lambda_n$
  denote the lowest nonzero eigenvalue of the Laplacian on~$E_n$.
  Assume that there exists a constant $c>0$ such that $\lambda_n \ge
  c$ for all~$n$.  Thus the sequence $\{E_n\}$ is an \emph{expanding
  sequence of graphs}.  The coarse space~$E$ provides a
  counter-example to the coarse Baum-Connes conjecture
  (see~\cite{HigsonSkandalisLafforgue}).  Let~$P$ be the spectral
  projection for the Laplacian, which has been shown is not in the
  range of the coarse Baum-Connes assembly map.  The above argument
  shows that the class of $[P]$ in $\K_0(C^*(E))$ pairs trivially with
  $\tilde\K_0(\bv(E))$.  More generally, if $X$ is a coarse space,
  $i\colon E \to X$ is a coarse embedding and~$\theta$ is any class in
  $\tilde{\K}_*\bigl(\bv(X)\bigr)$, then $\braket{i_*[P]}{\theta} =
  \braket{[P]}{i^*(\theta)} = 0$ by functoriality of the pairings and
  by our discussion above.  Hence such counter-examples cannot be
  detected by pairing with the stable Higson corona.  In particular,
  this discussion applies to the groups containing an expanding
  sequence of graphs constructed by Gromov in~\cite{Gro}.
\end{rmk}

\begin{ex}  \label{dualcoarsebaumconnesforwellspacedray}
  Consider now the special case $X=\bigsqcup X_n$ where each~$X_n$ is
  just a point (the \emph{well-spaced ray}).  Again the pairing
  between $\K_*\bigl(C^*(X)\bigr)$ and $\tilde\K_*\bigl(\bv(X)\bigr)$
  vanishes.  In this case, both the ordinary coarse assembly map~$\mu$
  and the coarse co\brd{}assembly map~$\mu^*$ are isomorphisms.  The
  assembly map is treated in~\cite{Yu2}, the co\brd{}assembly map can
  be treated easily using $\vv(X,D)=\ell^\infty(\N,D\otimes\Comp)$.
\end{ex}

\section{Homotopy invariance in the coefficient algebra}
\label{sec:functor_in_coefficient}

In order to prove that the coarse co\brd{}assembly map $\mu^*_{X,D}$
is an isomorphism for scalable spaces, we have to investigate the
homotopy invariance properties of $\tilde{\K}_*\bigl(\bv(X,D)\bigr)$
as a functor of~$D$.  For this, we have to assume the coarse structure
to be countably generated.  Since any countably generated coarse
structure arises from a metric as in
Example~\ref{themetriccoarsestructure}, we may fix such a metric.

A $*$\nbd{}homomorphism $f\colon D_1\to D_2$ induces compatible
$*$\nbd{}homomorphisms
\begin{align*}
  f\otimes 1 &\colon C_0(X,D_1\otimes\Comp) \to
  C_0(X,D_2\otimes\Comp),
  \\
  \bar{f}_X &\colon \vv(X,D_1) \to \vv(X,D_2),
  \\
  f_X &\colon \bv (X,D_1) \to \bv (X,D_2).
\end{align*}
Thus the assignments $D\mapsto \vv(X,D), \bv(X,D)$ are functors, which
we denote by $\vv(X,\cdot)$ and $\bv(X,\cdot)$, respectively.

\begin{prop}  \label{pro:stable_split_exact}
  The functors $\K\circ \vv(X,\cdot)$ and $\K\circ \bv(X,\cdot)$ are
  stable, split exact functors from the category of
  $C^*$\nbd{}algebras and $C^*$\nbd{}algebra homomorphisms to the
  category of Abelian groups and Abelian group homomorphisms,
  where~$\K$ denotes the $\K$\nbd{}theory functor.
\end{prop}

\begin{proof}
  If~$f$ is a completely bounded linear map $D_1\to D_2$, then it
  induces completely bounded linear maps $\bar{f}_X$ and~$f_X$ as
  above.  Hence our two functors map extensions of $C^*$\nbd{}algebras
  with completely bounded sections again to such extensions.  This
  yields the asserted split exactness.  Stability is already built
  into the definition.
\end{proof}

Lemma~\ref{injectivityofinclusions} implies easily that the functors
$\tilde{\K}\circ \bv(X,\cdot)$ and $\tilde{\K}\circ \vv(X,\cdot)$ are
also split exact and stable.

The main result of this section is that the functors $\K \circ \bv (X,
\cdot)$ and $\K \circ \vv (X, \cdot)$ are homotopy invariant.  That
is, if $f$ and~$f'$ are homotopic $*$\nbd{}homomorphisms $D_1\to D_2$,
then they induce the same maps $\K_*\bigl(\vv(X,D_1)\bigr)\to
\K_*\bigl(\vv(X,D_2)\bigr)$ and $\K_*\bigl(\bv(X,D_1)\bigr)\to
\K_*\bigl(\bv(X,D_2)\bigr)$.  This is not obvious because
$\vv(X,D\otimes C[0,1])$ and $\vv(X,D)\otimes C[0,1]$ are different.
Each of these involves functions on $X\times [0,1]$ with values
in~$D$, but the variation condition is slightly different.  This is
closely related to the difference between homotopy and operator
homotopy in Kasparov theory.  This becomes manifest in the particular
case when~$X$ is a discrete, countable group.  It is shown
in~\cite{EmersonMeyer} that for any $C^*$\nbd{}algebra~$D$ with
trivial $G$\nbd{}action, there is a canonical isomorphism
$$
\tilde{\K}_*\bigl(\bv(G,D)\bigr)
\cong \KK^G_{*+1}(\C,C_0(G)\otimes D).
$$
This is proved by showing that both groups are defined by the same
cycles and that the homotopy relation on $\K$\nbd{}theory corresponds
to operator homotopy in Kasparov's theory.  Thus the proof uses the
equivalence of operator homotopy and homotopy.  The homotopy
invariance of $\KK^G$ now implies the homotopy invariance of the
functor $\tilde\K\circ\bv(G,\cdot)$.  In order to extend this from
groups to general metric spaces, we have to adapt Kasparov's Technical
Theorem.  Before we do this, we note the following fact:

\begin{lemma}
  If the functor $\tilde\K\circ\bv(X,\cdot)$ is homotopy invariant,
  then so are $\K\circ\bv(X,\cdot)$, $\K\circ\vv(X,\cdot)$ and
  $\tilde\K\circ\vv(X,\cdot)$.
\end{lemma}

\begin{proof}
  Since $\K(\cdot)$ is homotopy invariant, the homotopy invariance of
  the reduced and unreduced theories is equivalent.  Since $\K\circ
  C_0(X,\cdot)$ is also homotopy invariant, homotopy invariance of
  $\K\circ\vv(X,\cdot)$ and $\K\circ\bv(X,\cdot)$ are equivalent by
  the Five Lemma.
\end{proof}

\begin{note}
  For the remainder of this article, we let $I\defeq C([0,1])$.
\end{note}

For a $C^*$\nbd{}algebra~$D$, we let $\ev_{D,t}\colon D\otimes I\to D$
be evaluation at $t\in I$.  A functor~$F$ is homotopy invariant if and
only if $F(\ev_{D,0})=F(\ev_{D,1})$ for all~$D$.  Since $\KK$ is
homotopy invariant, $\ev_0$ and $\ev_1$ define the same element in
$\KK_0(I,\C)$.  Hence the homotopy invariance results follow from the
following theorem.

\begin{thm}  \label{actionofkk}
  Let~$X$ be a countably generated coarse space and let~$D$ be a
  $C^*$\nbd{}algebra.  Then there is a well-defined pairing
  $$
  \tilde\K_*\bigl(\bv(X,D\otimes I)\bigr) \times \KK_0(I,\C) \to
  \tilde\K_*\bigl(\bv(X,D)\bigr),
  \qquad
  (x,y) \mapsto x \cdot y,
  $$
  which satisfies $x \cdot [\ev_t] = (\ev_{D,t})_*(x)$ for all
  $t\in[0,1]$.
\end{thm}

By Lemma~\ref{lem:discretization} and
Proposition~\ref{functorialityofcoronaalgebra}, we may assume~$X$ to
be discrete without loss of generality.  Hence we assume this in the
following.

A cycle~$b$ for $\KK_0(I,\C)$ is a pair $(\mathfrak{h}_2,F_2)$,
where~$\mathfrak{h}_2$ is a graded Hilbert space carrying a
$*$\nbd{}representation of~$I$ by even operators and where~$F_2$ is an
odd, self-adjoint operator on~$\mathfrak{h}_2$ for which $F^2_2-1$ and
$[f,F_2]$ for all $f\in I$ are compact.  We may assume that the
representation of~$I$ is essential, that is, $I\cdot
\mathfrak{h}_2=\mathfrak{h}_2$.

Proposition~\ref{identificationofreducedktheory} yields
$\tilde\K_*(\bv(X,\cdot))\cong \K_*(\bvr(X,\cdot))$.  We only consider
the pairing with $\K_*(\bvr(X,D\otimes I))$ for $*=1$, the case $*=0$
is similar.  Any class in $\K_1\bigl(\bvr(X,D\otimes I)\bigr)$ can be
represented by a unitary in $\bvr(X,D\otimes I)$, which we may then
lift to an element~$U$ of $\vvr(X,D\otimes I)$.  Thus $U\colon
X\to\Bound(D\otimes I\otimes H_1)$ is a bounded map of vanishing
variation such that $U(x)-U(y)\in \Comp(D\otimes I\otimes H_1)$ for
all $x,y\in X$, and such that $U^*U-1$ and $UU^*-1$ lie in
$C_0(X,\Comp(D\otimes I\otimes H_1))$.  Here~$H_1$ is another
(ungraded) Hilbert space with a countable basis.  Let
$$
\mathfrak{h}_1 = H_1 \oplus H_1^{\mathrm{op}},
\qquad
F_1 =
\begin{pmatrix} 0 & U^* \\ U & 0 \end{pmatrix}.
$$
Here $H_1 \oplus H_1^{\mathrm{op}}$ denots the Hilbert space
$H_1\oplus H_1$ graded by $\bigl(\begin{smallmatrix} 1 & 0 \\ 0 &
-1\end{smallmatrix}\bigr)$.  Thus~$F_1$ is a bounded map of vanishing
variation from~$X$ to odd, self-adjoint operators on the graded
Hilbert $D\otimes I$\nbd{}module $D\otimes I\otimes
\mathfrak{h}_1$, which also satisfies $F_1(x)-F_1(y)\in \Comp(D\otimes
I\otimes\mathfrak{h}_1)$ for all $x,y\in X$ and $F_1^2-1\in
C_0(X,\Comp (D\otimes I\otimes\mathfrak{h}_1))$.  The last property
mean that we have constructed a cycle for $\KK_0(\C,C_0(X)\otimes D
\otimes I)$.

The above construction produces from a unitary $a\in\bvr(X,D\otimes
I)$ a cycle~$\hat{a}$ for $\KK_0(\C,C_0(X)\otimes D \otimes I)$ with
the additional property that it has vanishing variation and is
constant up to compact perturbation as a function of~$X$.  Let us call
such cycles \emph{special}.  The only choice in this construction
of~$\hat{a}$ is the lifting of~$a$ to $U\in \vvr(X,D\otimes I)$.
Thus~$\hat{a}$ is determined uniquely up to operator homotopy.
Conversely, any special cycle~$\hat{a}$ for $\KK_0(\C,C_0(X,D))$ comes
from a unique unitary in $\bvr(X,D)$.  It is easy to see that the
equivalence relation of stable homotopy for unitaries in $\bvr(X,D)$
corresponds exactly to the equivalence relation generated by addition
of degenerate cycles and operator homotopy within the class of special
cycles for $\KK_0(\C,C_0(X,D))$.

Let~$\hat{a}$ be a special cycle for $\KK_0(\C,C_0(X,D\otimes I))$ as
above and let~$b$ be a cycle for $\KK_0(I,\C)$.  Our task is to show
that we can represent the Kasparov product $\hat{a}\otimes_I b\in
\KK_0(\C,C_0(X,D))$ again by a special cycle and that this special
representative is unique up to an operator homotopy among special
cycles.  Viewing the resulting special Kasparov cycles as elements of
$\K_1(\bvr(X,D))$, we obtain the desired pairing as in
Theorem~\ref{actionofkk} and hence the homotopy invariance result that
we are aiming for.

The referee has suggested to prove Theorem~\ref{actionofkk} using the
groupoid $G(X)$ associated to the coarse space~$X$
in~\cite{SkandalisYuTu}.  Special cycles for $\KK(\C,C_0(X,D))$ are
the same as cycles for $\KK^{G(X)}(\ell^\infty(X), C_0(X,D))$.
However, since neither $\ell^\infty(X)$ nor $G(X)$ are separable and
since $\ell^\infty(X)$ occurs in the first variable, some work is
still needed to get Kasparov products in this situation.  We follow a
different route in the following.

Let $F_1$ and $F_2$, $\mathfrak{h}_1$ and $\mathfrak{h}_2$ be as in
the discussion above.  Let $F_1\sharp F_2$ be the collection of
bounded maps $F\colon X\to
\Bound(D\otimes\mathfrak{h}_1\hat{\otimes}\mathfrak{h}_2)$ taking
values in odd, self-adjoint operators and satisfying the following
conditions:
\begin{enumerate}[(1)]
\item $F(x)-F(y)\in \Comp(D\otimes\mathfrak{h}_1 \hat{\otimes}
  \mathfrak{h}_2)$ for all $x,y\in X$.

\item $F$ has vanishing variation.  

\item $F^2-1\in C_0(X,\Comp(D\otimes\mathfrak{h}_1 \hat{\otimes}
  \mathfrak{h}_2))$.

\item $F$ is an $F_2$ connection.  

\item $[F_1\hat{\otimes} 1,F]\ge 0$ modulo
  $C_0(X,\Comp(D\otimes\mathfrak{h}_1 \hat{\otimes} \mathfrak{h}_2))$.

\end{enumerate}

Conditions (3)--(5) say that~$F$ is a Kasparov product of $F_1$
and~$F_2$ in the usual sense and (1)--(2) say that~$F$ is special.
Thus we have to show that $F_1\sharp F_2$ is not empty and that any
two elements of $F_1\sharp F_2$ are operator homotopic.

\begin{prop}[A coarse Technical Theorem]
  \label{coarsetechnicaltheorem}
  The set $F_1\sharp F_2$ is non-empty.  
\end{prop}

\begin{proof}
  We follow Kasparov's argument in~\cite{Kasparov}.  We let
  \begin{align*}
    A_1 &\defeq C_0(X,\Comp (D\otimes\mathfrak{h}_1)\hat{\otimes} 1)
    + C_0(X,\Comp (D\otimes\mathfrak{h}_1\hat{\otimes} \mathfrak{h}_2)),
    \\
    A_2 &\defeq C^*\bigl(1\hat{\otimes} (F^2_2-1)\bigr)
    \subseteq 1\otimes \Comp(\mathfrak{h}_2),
    \\
    J &\defeq
    C_0(X,\Comp (D\otimes\mathfrak{h}_1 \hat{\otimes} \mathfrak{h}_2)).
  \end{align*}
  These are separable $C^*$\nbd{}subalgebras of
  $\Bound(D\otimes\mathfrak{h}_1 \hat{\otimes} \mathfrak{h}_2)$.  Let
  $h_1$, $h_2$ and~$k$ be strictly positive elements of $A_1$, $A_2$
  and~$J$, respectively.

  Let $\Delta_1\subseteq\Bound(\mathfrak{h}_1)$ be the closed linear
  span of the operators $F_1(x)$, $x\in X$.  Let $\Delta\subseteq
  \Bound(D\otimes\mathfrak{h}_1 \hat{\otimes}\mathfrak{h}_2)$ be the
  closed linear span of
  $$
  \{F_1(x)\otimes 1 \mid x\in X\} \cup\{1\otimes F_2\} \cup
  \{h_1(x) \mid x\in X\} \cup \{h_2\}.
  $$
  These spaces are separable.  Hence there exist compact subsets
  $Y_1\subseteq \Delta_1$ and $Y\subseteq\Delta$ such that $\C\cdot
  Y_1$ and $\C\cdot Y$ are dense in $\Delta_1$ and~$\Delta$,
  respectively.  We may assume that $h_1(x)\in Y$ for all $x\in X$
  because $h_1(X)$ is compact.

\begin{lemma}  \label{existenceofacertainapproximateunit}
  There exists an approximate unit $(e_i)_{i\in\N}$ for the
  $C^*$\nbd{}algebra $\Comp(D\otimes\mathfrak{h}_1)$ with $0\le e_i
  \le e_{i+1} \le 1$ for all~$i$ such that all~$e_i$ are even and
  $\lim_{i\to\infty} \norm{[e_i,y]}=0$ uniformly for $y\in Y_1$.

  There exists an approximate unit $(\epsilon_i)_{i\in\N}$ for the
  $C^*$\nbd{}algebra $\Comp (D \otimes \mathfrak{h}_1 \hat{\otimes}
  \mathfrak{h}_2)$ with $0 \le \epsilon_i \le \epsilon_{i+1} \le 1$
  for all~$i$ such that all~$\epsilon_i$ are even and
  $\lim_{i\to\infty} \norm{[\epsilon_i,y]}=0$ uniformly for $y\in Y$.
\end{lemma}

\begin{proof}
  Both assertions are special cases of \cite{Kasparov}*{Lemma 1.4}.
  For the first one, apply Kasparov's lemma with the trivial
  group~$G$, $A\defeq\Comp(D\otimes\mathfrak{h}_1)$,
  $B\defeq\Bound(D\otimes\mathfrak{h}_1)$, $Y\defeq Y_1$, and the
  inclusion map $\varphi\colon Y\to B$.
\end{proof}

Now let~$(\rho_i)$ be a sequence of functions in $C^1_c(\R_+)$ that
satisfy
\begin{enumerate}[(1)]
\item $0\le\rho_i \le\rho_{i+1}\le 1$;

\item $\rho_i(t)=1$ for $t\le i$;

\item $\abs{\rho_i'(t)}\le 2^{-i}$ for all $t\in[0,\infty)$.

\end{enumerate}
The last condition implies $\Var_N \rho_i(t) \le 2^{-i}N$ for all
$t\in\R_+$, $i,N\in\N$.  Let~$x_0$ be a point of~$X$ and define
functions~$\psi_i$ on~$X$ by
$$
\psi_i(x)\defeq \rho_i(d(x,x_0)).
$$
The triangle inequality yields $\Var_N \psi_i(x) \le 2^{-i}N$ for all
$x\in X$, $i,N\in\N$.  Moreover, $\psi_i$ equals~$1$ on $B_i(x_0)$, is
compactly supported and satisfies $0\le \psi_i\le 1$.  Define
$$
u_i (x) \defeq \psi_i (x) e_i,
\qquad
v_i(x) \defeq \psi_i (x) \epsilon_i,
$$
where $(e_i)$ and $(\epsilon_i)$ are as in
Lemma~\ref{existenceofacertainapproximateunit}.  Thus $(u_i)$ is an
approximate unit for~$A_1$ and $(v_i)$ is an approximate unit
for~$J$.  Moreover, we have
$$
\Var_N u_i(x) \le 2^{-i}N,
\qquad
\Var_N v_i(x) \le 2^{-i}N
\qquad \text{for all $N, i\in\N$, $x\in X$.}
$$

Passing to a subsequence, we can achieve
\begin{displaymath}
  \norm{u_i h_1 - h_1} \le 2^{-i}
  \qquad \text{for all $i\in\N$}
\end{displaymath}
because $h_1\in A_1$.  For $y\in Y_1\subseteq
\Bound(\mathfrak{h}_1)\subseteq
\Bound(D\otimes\mathfrak{h}_1\hat{\otimes}\mathfrak{h}_2)$ and $x\in
X$, the commutator $[u_i(x),y]$ is given by the function
$$
X\to \Bound(D\otimes\mathfrak{h}_1\hat{\otimes} \mathfrak{h}_2),
\qquad x\mapsto [u_i(x),y] = \psi_i(x)[e_i,y],
$$
and $\norm{[e_i,y]}\to0$ uniformly on~$Y_1$ as $i\to\infty$.
Therefore, passing again to a subsequence we can achieve
\begin{displaymath}
  \norm{[u_i(x),y]} \le 2^{-i},
  \qquad \text{for all $x\in X$, $y\in Y_1$, $i\in\N$.}
\end{displaymath}
We still have
\begin{displaymath}
  \Var_N u_i (x) \le 2^{-i}N,
  \qquad \text{for all $x\in X$, $i,N\in\N$.}
\end{displaymath}

Since $(v_i)$ is an approximate unit for~$J$ and $u_ih_2\in J$, we
have $\norm{v_i u_ih_2 - u_ih_2} \to 0$ as $i\to\infty$.  Passing to a
subsequence of $(v_i)$, we can achieve
\begin{displaymath}
  \norm{v_iu_ih_2 - u_ih_2} \le 2^{-2i},
  \qquad
  \norm{v_i u_{i+1}h_2 - u_{i+1}h_2} \le 2^{-2i},
  \qquad
  \norm{v_i k - k} \le 2^{-2i}.
\end{displaymath}

Now recall that $\norm{[v_i(x),y]} \le \norm{[\epsilon_i,y]}\to 0$
uniformly for $x\in X$, $y\in Y$ as $i\to\infty$.  Regarding elements
of~$Y$ as constant functions on~$X$, we obtain $\norm{[v_i,y]}\to 0$
uniformly for $y\in Y$ as $i\to\infty$.  Since $h_1(X)\subseteq Y$, we
also obtain $\norm{[v_i,h_1]}\to 0$ as $i\to\infty$.  Let $b_i\defeq
(v_i-v_{i-1})^{1/2}$.  Passing once again to a subsequence, we can
achieve
\begin{displaymath}
  \norm{[b_i,h_1]} \le 2^{-i},
  \qquad
  \norm{[b_i,h_2]} \le 2^{-i},
  \qquad
  \norm{[v_i,y]} \le 2^{-i},
  \qquad
  \norm{[b_i,y]} \le 2^{-i},
\end{displaymath} 
for all $y\in Y$.  Finally, we can assume
\begin{displaymath}
  \Var_N v_i(x) \le 2^{-i}N,
  \qquad
  \Var_N b_i(x) \le 2^{-i}N
  \qquad \text{for all $i,N\in\N$, $x\in X$.}
\end{displaymath}

Since~$b_i$ is a compactly supported function
$X\to\Comp(D\otimes\mathfrak{h}_1\hat{\otimes}\mathfrak{h}_2)$, so is
$q_i\defeq b_iu_ib_i$.  Moreover, the $N$\nbd{}variation of~$q_i$ is
at most $3\cdot 2^{-i}N$ everywhere on~$X$.

As in~\cite{Kasparov}, the equations established above imply that the
series $\sum_{i=1}^{\infty} q_i$ converges in the strict topology to
some $M_2\in\Mult(J)$.  Equivalently, $\sum_{i=1}^\infty q_i a$ is
norm convergent for every $a\in J$.  Let $M_1\defeq 1-M_2$.  One next
checks that the series $\sum_{i=1}^\infty q_ih_2$ is also absolutely
convergent, so that $M_2 h_2\in J$.  Similarly, $M_1 h_1 \in J$ and
$[M_2,y]\in J$ for all $y\in Y$, so that also $[M_1,y] \in J$ for all
$y\in Y$.  We claim that $M_1$ and~$M_2$ have vanishing variation.  To
check this, choose $\epsilon>0$ and $R\ge 0$.  Choose~$m$ large enough
that $\sum_{i=m}^\infty 2^{-i}R <\epsilon/6$.  The function
$\sum_{i=1}^{m-1} \norm{q_i}$ is a $C_0$\nbd{}function.  Hence we may
choose a compact set $K\subseteq X$ such that if $z\in X\setminus K$
then $\sum_{i=1}^{m-1} \norm{q_i(z)} < \epsilon/4$.  If $d(x,K)>R$ and
$d(x,y)\le R$, then
\begin{multline*}
  \norm{M_2(x)-M_2(y)}
  \le \sum_{i=1}^\infty \norm{q_i(x)-q_i(y)}
  = \sum_{i=m}^\infty \norm{q_i(x)-q_i(y)}
  + \sum_{i=1}^{m-1} \norm{q_i(x)-q_i(y)}
  \\ \le \sum_{i=m}^\infty 3\cdot 2^{-i}R
  + \sum_{i=1}^{m-1} \norm{q_i(x)} + \norm{q_i(y)}
  < \epsilon/2 + \epsilon/4 + \epsilon/4
  < \epsilon.
\end{multline*}
That is, $M_2$ has vanishing variation.  A similar estimate shows that
the series $M_2(x)-M_2(y) = \sum_{i=1}^\infty q_i(x)-q_i(y)$ converges
absolutely.  Since the summands are compact, so is $M_2(x)-M_2(y)$.
That is, $M_2\in \vvr(X,D)$.  Therefore,
$$
F\defeq M_1^{1/2}(F_1\otimes 1) + M_2^{1/2}(1\otimes F_2)
$$
also belongs to $\vvr(X,D)$.  The remaining conditions (3)--(5) for
$F\in F_1\sharp F_2$ are checked as in~\cite{Kasparov}.
\end{proof}

\begin{lemma}
  Any two elements of $F_1 \sharp F_2$ are operator homotopic.
\end{lemma}

\begin{proof} \emph{(Sketch)}
  Let $F$ and~$F' $ be two elements of $F_1 \sharp F_2$. Following the
  method of the previous proof, one constructs bounded functions
  $N_i\colon X\to \Bound(D\otimes\mathfrak{h}_1\hat{\otimes}
  \mathfrak{h}_2)$ of vanishing variation with $N_1+N_2=1$ such that
  $F'' \defeq N_1^{1/2} (F_1\hat{\otimes}1) + N_2^{1/2} F$ lies in
  $F_1 \sharp F_2$ and $[F,F''] \ge 0$ and $[F',F''] \ge 0$ modulo the
  ideal $C_0(X,\Comp(D\otimes \mathfrak{h}_1\hat{\otimes}
  \mathfrak{h}_2)$.  One now follows a standard argument
  (\cite{Bla}*{p.\ 149}) and writes $[F,F''] = P + K$, where $P$ is a
  positive, even element of $\vvr(X,\Comp(D\otimes
  \mathfrak{h}_1\hat{\otimes} \mathfrak{h}_2))$ and where $K\in
  C_0(X,\Comp(D\otimes \mathfrak{h}_1\hat{\otimes} \mathfrak{h}_2))$.
  Now the path $G_t\defeq (1+\cos(t)\sin(t)P)^{-1/2}(\cos(t)F +
  \sin(t)F'')$ is an operator homotopy between $F$ and~$F''$ in
  $F_1\sharp F_2$.  A similar formula produces an operator homotopy
  between $F'$ and~$F''$.  Hence $F$ and~$F'$ are operator homotopic
  in $F_1\sharp F_2$ as desired.
\end{proof}

Next once has to check the following facts:
\begin{enumerate}[(1)]
\item if there exists an operator homotopy between $F_1$ and~$F_1'$
  consisting of special cycles for $\KK(\C,C_0(X,D\otimes I))$, then
  there is an operator homotopy between appropriate elements of
  $F_1\sharp F_2$ and $F_1'\sharp F_2$ consisting of special cycles
  for $\KK(\C,C_0(X,D))$;
  
\item if $F_2$ and~$F_2'$ are operator homotopic, then there is an
  operator homotopy between appropriate elements of $F_1\sharp F_2$
  and $F_1\sharp F_2'$ consisting of special cycles for
  $\KK(\C,C_0(X,D))$;

\item if $F_1$ or $F_2$ is degenerate, then $F_1\sharp F_2$ contains a
  degenerate cycle.

\end{enumerate}
The homotopies required for the first two assertions are again
constructed as Kasparov products, and the third assertion is
well-known.  This establishes that the pairing that sends $F_1,F_2$ to
any element of $F_1\sharp F_2$ descends to operator homotopy
equivalence classes.

This finishes the proof of Theorem~\ref{actionofkk} and hence of the
homotopy invariance of the functor $\tilde\K\circ \bv(X,\cdot)$.  As
explained, this implies the homotopy invariance of the functors
$\tilde\K\circ\vv(X,\cdot)$, $\K\circ\bv(X,\cdot)$ and
$\K\circ\vv(X,\cdot)$.  By the universal properties of $\KK$, it
follows that an element of $\KK(D_1,D_2)$ induces canonical maps
$F(D_1)\to F(D_2)$, where~$F$ is any one of these functors.  To
summarize, we have:

\begin{cor}  \label{cor:KK_factorization}
  The functors $\tilde\K\circ\bv(X,\cdot)$,
  $\tilde\K\circ\vv(X,\cdot)$, $\K\circ\bv(X,\cdot)$ and
  $\K\circ\vv(X,\cdot)$ factor through the category $\KK$.
\end{cor}

\section{Maps of weakly vanishing variation and scalable spaces}
\label{sec:scalable}

In this section we apply the above ideas to prove that the map
$\mu^*_{X,D}$ is an isomorphism for scalable, uniformly contractible
metric spaces.  In fact we prove that $\tilde{\K}_*
\bigl(\vv(X,D)\bigr) = 0$ for any scalable space~$X$ and any
$C^*$\nbd{}algebra~$D$.  We begin by discussing a type of
functoriality of the algebras $\vv(X,D)$ in the $X$\nbd{}variable.

\begin{defn}  \label{definitionofwvvmap}
  Let~$X$ be a coarse space, let~$Y$ be a coarse metric space, and let
  $f\colon X \to Y$ be a continuous map.  We say that~$f$ has
  \emph{weakly vanishing variation}, or is WVV, if~$f$ maps entourages
  to entourages and if
  $$
  f^{-1}(K) \cap (\Var_E f)^{-1}([\epsilon,\infty))
  $$
  is bounded for all $\epsilon>0$, all entourages $E\subseteq X\times
  X$ and all bounded subsets $K\subseteq Y$.
\end{defn}

This allows us to treat coarse maps and maps of vanishing variation
simultaneously because of the following obvious lemma.

\begin{lemma}
  Let $X$, $Y$ and~$f$ be as in Definition~\ref{definitionofwvvmap}.
  If~$f$ is coarse or if~$f$ has vanishing variation, then~$f$ has
  weakly vanishing variation.
\end{lemma}

\begin{lemma}
  Let $X$ be a coarse space, $Y$ a coarse metric space, and~$Z$ a
  metric space, and let $\phi\colon X\to Y$ and $\psi\colon Y\to Z$ be
  continuous maps.  If~$\psi$ has vanishing variation and~$\phi$ has
  weakly vanishing variation, then $\psi\circ\phi\colon X\to Z$ has
  vanishing variation.
\end{lemma}

\begin{proof}
  Let $\epsilon>0$ and let~$E$ be an entourage in~$X$.  We must show
  that there exists a bounded subset $K\subseteq X$ such that
  $\Var_E(\psi\circ\phi)(x)<\epsilon$ for $x\in X\setminus K$.  Since
  $\phi(E)$ is an entourage in~$Y$, there is $R>0$ with
  $d(\phi(x),\phi(x'))<R$ for $(x,x')\in E$.  Choose a bounded subset
  $L\subseteq Y$ such that $\Var_R\psi(y) <\epsilon$ for $y\in
  Y\setminus L$.  Since~$\psi$ has vanishing variation and is
  continuous, it is uniformly continuous.  Hence there is $\delta>0$
  such $d(\psi(y),\psi(y')) < \epsilon$ whenever $d(y,y')<\delta$.
  Now let $K\defeq\phi^{-1}(L)\cap(\Var_E\phi)^{-1}([\delta,\infty))$.
  Then~$K$ is bounded since~$\phi$ has weakly vanishing variation.
  Choose $(x,x')\in E$ with $x\in X\setminus K$.  Then $\phi(x)\in
  Y\setminus L$ or $\Var_E \phi(x)<\delta$.  Suppose that $\phi(x)\in
  Y\setminus L$.  Since also $d(\phi(x),\phi(x'))<R$, we get
  $d(\psi\phi(x),\psi\phi(x')) < \epsilon$ by choice of~$L$.  If
  $\Var_E \phi(x)<\delta$, then $d(\phi(x),\phi(x'))<\delta$ and thus
  $d(\psi\phi(x),\psi\phi(x'))<\epsilon$.  That is, $\psi\phi$ has
  vanishing variation.
\end{proof}

\begin{cor}  \label{newfunctoriality}
  Let $X$ be a coarse space and let~$Y$ be a coarse metric space.
  Then a continuous map $\phi\colon X\to Y$ of weakly vanishing
  variation induces $*$\nbd{}homomorphisms $\vv(Y,D)\to\vv(X,D)$ and
  $\vvr(Y,D)\to\vvr(X,D)$ by the formula $f\mapsto f\circ\phi$.
\end{cor}

Hence $\K(\vv({\cdot},D))$ and $\tilde\K(\vv({\cdot},D))$ are
functorial for WVV maps.  However, since WVV maps need not be proper,
they need \emph{not} act on $\K(C_0({\cdot},D))$, $\K(\bv({\cdot},D))$
or $\K(\bvr({\cdot},D))$.

Let $X$ be a coarse space, let~$Y$ be a coarse metric space, and let
$\phi,\phi'\colon X\to Y$ be WVV maps.  Let $\Phi\colon
X\times[0,1]\to Y$ be a continuous map with 
$\Phi(x,0) = \phi(x)$ and $\Phi(x,1) = \phi'(x)$ for all $x\in X$.
For an entourage~$E$ in~$X$, we define
$$
\Var_E^1 \Phi \colon X\times[0,1] \to [0,\infty),
\qquad
(x,t) \mapsto \sup\{d(\Phi(x,t),\Phi(y,t)) \mid (x,y)\in E\}.
$$

\begin{defn}  \label{def:WVV_homotopy}
  We say that~$\Phi\colon X\times[0,1]\to Y$ is a \emph{WVV homotopy}
  between $\phi$ and~$\phi'$ and call $\phi$ and~$\phi'$ \emph{WVV
    homotopic} if the function $\Var_E^1$ is bounded for any
  entourage~$E$ and the set $\Phi^{-1}(K)\cap(\Var_E^1
  \Phi)^{-1}([\epsilon,\infty))$ is bounded in $X\times[0,1]$ for all
  $\epsilon>0$, all entourages $E\subseteq X\times X$ and all bounded
  subsets $K\subseteq Y$.  We call~$Y$ \emph{WVV contractible} if the
  identity on~$Y$ is WVV homotopic to a constant map.
\end{defn}

\begin{prop}  \label{mainconsequenceofwvvhomotopy}
  Let $X$ and~$Y$ be as above, let $\Phi\colon X\times[0,1]\to Y$ be a
  WVV homotopy and let~$D$ be a $C^*$\nbd{}algebra.  Then $f\mapsto
  f\circ\Phi$ defines $*$\nbd{}homomorphisms $\vv(Y,D)\to
  \vv(X,D\otimes I)$ and $\vvr(Y,D)\to \vvr(X,D\otimes I)$.
\end{prop}

\begin{proof}
  Continuity of $\Phi$ and~$f$ imply immediately that $t\mapsto
  f(\Phi(x,t))$ indeed lies in $D\otimes I = C([0,1],D)$ for every
  $x\in X$ and that the map $f\circ\Phi\colon X\to D\otimes I$ is
  continuous.  It is evidently bounded.  It remains to check the
  variation condition.  Choose an entourage~$E$ in~$X$ and
  $\epsilon>0$.  We have to find a bounded subset $K\subseteq X$ such
  that $\norm{f(\Phi(x,t))-f(\Phi(x',t))}<\epsilon$ for $(x,x')\in E$,
  $t\in[0,1]$ with $x\notin K$.
  
  Since continuous vanishing variation functions are uniformly
  continuous, we can find $\delta>0$ such that
  $\norm{f(y)-f(y')}<\epsilon$ if $d(y,y')<\delta$.  Since $\Var_E^1$
  is bounded, there is $R>0$ such that $d(\Phi(x,t),\Phi(x',t))\le R$
  for all $(x,x')\in E$, $t\in[0,1]$.  Since~$f$ has vanishing
  variation, we can find a bounded subset $L\subseteq Y$ such that
  $\Var_R f(y)<\epsilon$ for $y\in Y\setminus L$.  Let $K\subseteq X$
  denote the projection to~$X$ of the bounded subset
  $$
  \Phi^{-1}(L) \cap(\Var_E^1 \Phi)^{-1}([\delta,\infty))
  \subset X\times[0,1].
  $$
  If $x\in X \setminus K$, $t\in[0,1]$, then $\Phi(x,t)\in Y\setminus
  L$ or $\Var_E^1 \Phi(x,t)<\delta$.  Suppose first that $\Phi(x,t)
  \in Y\setminus L$.  Since $(x,x')\in E$, we have
  $d(\Phi(x,t),\Phi(x',t))\le R$.  The choice of~$L$ yields
  $\norm{f(\Phi(x,t)) - f(\Phi(x',t))} < \epsilon$ as desired.  If
  $\Var_E^1 \Phi(x,t) < \delta$, then $d(\Phi(x,t),\Phi(x',t)) <
  \delta$ and hence $\norm{f(\Phi(x,t)) - f(\Phi(x',t))} < \epsilon$
  as well by the choice of~$\delta$.
\end{proof}
 
\begin{rmk}
  The notion of WVV homotopy is motivated by the following example.
  Let~$X$ be a complete Riemannian manifold of non\brd{}positive
  curvature.  Fix a point $x_0\in X$ and let $\exp\colon T_{x_0}(X)
  \to X$ be the exponential map at~$x_0$.  It is well known that
  $\exp$ is a diffeomorphism satisfying $d(\exp v,\exp w) \ge d(v,w)$
  for all tangent vectors $v,w \in T_{x_0}(X)$.  Let $\log$ denote the
  inverse of $\exp$.  Define
  $$
  \Phi\colon X \times[0,1] \to X,
  \qquad
  \Phi(x,t)\defeq \exp(t\log x).
  $$
  Then it is easy to check that~$\Phi$ is a WVV homotopy between
  the identity map $X\to X$ and the constant map~$x_0$.  That is, $X$
  is WVV contractible.
\end{rmk}

More generally, we can make the following definition.

\begin{defn}
  Let~$X$ be a coarse metric space.  We call~$X$ \emph{scalable} if
  there is a continuous map $r\colon X\times [0,1]\to X$ such that
  \begin{enumerate}[(1)]
  \item $r_1(x)=x$;
  \item the map $X\times [\epsilon,1]\to X$, $(x,t)\mapsto r_t(x)$ is
    proper for all $\epsilon>0$;
  \item the maps $r_t$ are uniformly Lipschitz and satisfy
    $$
    \lim_{t\to 0} \sup_{x\neq x'\in X}
    \frac{d(r_t(x),r_t(x'))}{d(x,x')} = 0.
    $$
  \end{enumerate}

\end{defn}

The following proposition follows immediately from the definition.

\begin{prop}  \label{scalableimplieswvvnullhomotopic}
  Scalable spaces are WVV contractible.
\end{prop}

\begin{cor}
  If~$X$ is a scalable space, then
  $\tilde{\K}_*\bigl(\vv(X,D)\bigr)=0$ for every
  $C^*$\nbd{}algebra~$D$.  If~$X$ is also uniformly contractible and
  has bounded geometry, then~$\mu^*$ is an isomorphism.
\end{cor}

\begin{proof}
  Since~$X$ is scalable, the identity map and a constant map on~$X$
  are WVV homotopic.  Proposition~\ref{mainconsequenceofwvvhomotopy}
  and Corollary~\ref{cor:KK_factorization} yield that they induce the
  same map on $\tilde\K_*(\vv(X,D))$.  Since we divided out the
  contribution of the constant functions in the reduced
  $\K$\nbd{}theory, we get $\tilde\K_*(\vv(X,D))\cong0$.  By
  Theorem~\ref{identificationofcoarsektheory} and the $\K$\nbd{}theory
  long exact sequence, this is equivalent to~$\mu^*$ being an
  isomorphism.
\end{proof}

We remark that Higson and Roe show in~\cite{HigsonRoe} that the coarse
Baum-Connes conjecture is an isomorphism for scalable spaces.

\begin{cor}
  Let~$\pi$ be the fundamental group of a compact, aspherical manifold
  of non\brd{}positive curvature.  Then the coarse co\brd{}assembly map
  $\mu^*_{\pi,D} \colon \tilde\K_{*+1}(\bv(\pi,D)\bigr) \to
  \KX^*(\pi,D)$ is an isomorphism for every coefficient
  $C^*$\nbd{}algebra~$D$.
\end{cor}

Similar results hold for groups~$G$ that admit cocompact, isometric,
proper actions on $\mathrm{CAT}(0)$ spaces.

\section{Groups which uniformly embed in Hilbert space}
\label{sec:embed_Hilbert}

We have introduced the co\brd{}assembly map in~\cite{EmersonMeyer}
because of its close relationship to the existence of a dual Dirac
morphism in the group case.  In this section, we use the notation of
\cites{EmersonMeyer,MeyerNest} concerning Dirac morphisms, dual Dirac
morphisms and $\gamma$\nbd{}elements.  We give a few explanations in
the proof of Proposition~\ref{pro:gamma}.

\begin{thm}[\cite{EmersonMeyer}]
  If a discrete group~$G$ has a dual Dirac morphism, then
  $$
  \mu^*_{G,D}\colon \tilde\K_{*+1}\bigl(\bv(G,D)\bigr)\to\KX^*(G,D)
  $$
  is an isomorphism for every $C^*$\nbd{}algebra~$D$.
\end{thm}

We are going to use this fact to prove that the coarse-coassembly map
is an isomorphism for groups that embed uniformly in a Hilbert space.
This method only applies to coarse spaces that are quasi-isometric to
a group.  It would be nice to have a more direct proof that~$\mu^*$ is
an isomorphism that applies to all coarse spaces that uniformly embed.

\begin{thm}  \label{uniformembeddingimpliesexistenceofgammaelement}
  Let~$G$ be a countable discrete group that embeds uniformly in a
  Hilbert space.  Then~$G$ possesses a dual Dirac morphism.  Hence the
  coarse co-assembly map for~$G$ is an isomorphism.
\end{thm}

Yu has shown the analogous result for the coarse assembly map for all
coarse spaces with bounded geometry (\cite{Yu:embeddable}).

\begin{rmk}
  It has been pointed out to us that Georges Skandalis and Jean-Louis
  Tu are aware of
  Theorem~\ref{uniformembeddingimpliesexistenceofgammaelement}; their
  work is independent of ours.
\end{rmk}

By a theorem of Higson, Guentner and Weinberger
(\cite{HigsonGuentnerWeinberger}), every countable subgroup of either
$\mathrm{GL}_n(k)$ for some field~$k$ or of an almost connected Lie
group admits a uniform embedding in Hilbert space.  Consequently one
obtains the following extension of Kasparov's results
in~\cite{Kasparov}:

\begin{cor}  \label{cor:gamma_for_linear}
  If~$G$ is a countable subgroup either of $\mathrm{GL}_n(k)$ for some
  field~$k$ or of an almost connected Lie group, then~$G$ possesses a
  dual Dirac morphism.
\end{cor}

The proof of
Theorem~\ref{uniformembeddingimpliesexistenceofgammaelement} is a
consequence of various results of Higson, Skandalis, Tu, and Yu (see
\cites{Hig2, SkandalisYuTu}).  Higson shows in~\cite{Hig2} that if~$G$
is a discrete group admitting a topologically amenable action on a
compact metrizable space, then the Novikov conjecture holds for~$G$.
But the argument manifestly also applies to the potentially larger
class of groups admitting an a\nbd{}T\brd{}menable action on a compact
space.

\begin{defn}[see~\cite{SkandalisYuTu}]  \label{def:negtype}
  Let~$G$ be a discrete group and~$X$ a compact $G$\nbd{}space.  We
  call the action of~$G$ on~$X$ \emph{a\nbd{}T\brd{}menable} if there
  exists a proper, continuous, real-valued function~$\psi$ on $X\times
  G$ satisfying:
  \begin{enumerate}[\ref{def:negtype}.1.]
  \item $\psi(x,e)=0$ for all $x\in X$;

  \item $\psi(x,g)=\psi(g^{-1}x,g^{-1})$ for all $x\in X$, $g\in G$;

  \item $\sum_{i,j=1}^n t_i t_j \psi(g^{-1}_i x,g^{-1}_i g_j)\le0$
    for all $x\in X$, $t_1,\ldots t_n\in\R$, $g_1,\ldots g_n\in G$
    for which $\sum t_i = 0$.

  \end{enumerate}
  Such a function~$\psi$ is called a \emph{negative type function}.
\end{defn}

The existence of a negative type function implies that the groupoid
$G\cross X$ possesses an affine isometric action on a continuous
family of Hilbert spaces over~$X$.  The above definition is relevant
because of the following results of \cites{SkandalisYuTu,Tu}:

\begin{thm}[\cite{SkandalisYuTu}]  \label{embeddable_acts}
  Let~$G$ be a discrete group.  If~$G$ admits a uniform embedding in a
  Hilbert space, then~$G$ admits an a\nbd{}T\brd{}menable action on a
  second countable compact space.
\end{thm}

\begin{thm}[\cite{Tu}]  \label{tustheorem}
  Let~$G$ be a discrete group and let~$X$ be a locally compact
  $G$\nbd{}space.  If~$G$ acts a\nbd{}T\brd{}menably on~$X$, then the
  transformation groupoid $G\cross X$ has a dual Dirac morphism and we
  have $\gamma=1$.
\end{thm}

We are not going to need the fact that $\gamma=1$.

If~$X$ is any compact space, let $\mathrm{Prob}(X)$ denote the
collection of probability measures on~$X$, equipped with the
weak-$^*$\brd{}topology.  This is again a compact space.  The space
$\mathrm{Prob}(X)$ is convex and hence contractible.  Even more, it is
equivariantly contractible with respect to any action of a compact
group on $\mathrm{Prob}(X)$.

\begin{lemma}[\cite{SkandalisYuTu}]  \label{inducedactionisamenable}
  Let~$X$ be a compact, second countable space on which a discrete
  group~$G$ acts a\nbd{}T\brd{}menably.  Then~$G$ acts
  a\nbd{}T\brd{}menably on $\mathrm{Prob}(X)$ as well.
\end{lemma}

\begin{prop}  \label{pro:gamma}
  Let~$G$ be a second countable locally compact group and let~$X$ be a
  second countable compact $G$\nbd{}space.  Suppose that~$X$ is
  $H$\nbd{}equivariantly contractible for all compact subgroups
  $H\subseteq G$ and that the groupoid $G\cross X$ has a dual Dirac
  morphism.  Then~$G$ has a dual Dirac morphism as well.
\end{prop}

\begin{proof}
  Let $\Dirac\in\KK^{G\cross X}(\Proj,C(X))$ be a \emph{Dirac
  morphism} for $G\cross X$ in the sense of~\cite{MeyerNest}.  This
  means two things: first, $\Dirac$ is a \emph{weak equivalence}, that
  is, for any compact subgroup $H\subseteq G$, restriction to~$H$
  maps~$\Dirac$ to an invertible morphism in $\KK^{H\cross
  X}(\Proj,C(X))$.  Secondly, $\Proj$ belongs to the localizing
  subcategory of $\KK^{G\cross X}$ that is generated by compactly
  induced algebras.  This subcategory contains all proper $G\cross
  X$-$C^*$\brd{}algebras.  It is shown in~\cite{MeyerNest} that a
  Dirac morphism for $G\cross X$ always exists.  A \emph{dual Dirac
  morphism} is an element $\eta\in \KK^{G\cross X}(C(X),\Proj)$ such
  that $\eta\circ\Dirac=\mathrm{id}_\Proj$.  It is shown in
  \cite{MeyerNest}*{Theorem 8.2} that a dual Dirac morphism exists
  whenever the Dirac dual Dirac method in the sense of proper actions
  applies.  In particular, it exists in the situation of
  Theorem~\ref{tustheorem}.

  Since~$X$ is compact, $C(X)$ contains the constant functions.  This
  defines a $G$\nbd{}equivariant $*$\nbd{}homomorphism $j\colon \C\to
  C(X)$.  The contractibility hypothesis on~$X$ insures that~$j$ is
  invertible in $\KK^H$ for any compact subgroup~$H$.  That is, $j$ is
  a weak equivalence.  Proposition 4.4 in~\cite{MeyerNest} yields
  that~$j$ induces an isomorphism
  $\KK^G(\Proj,\C)\cong\KK^G(\Proj,C(X))$.  Thus we obtain
  $\Dirac'\in\KK^G(\Proj,\C)$ with $j_*(\Dirac')=F(\Dirac)$, where
  $F\colon \KK^{G\cross X}\to \KK^G$ is the functor that forgets the
  $X$\nbd{}structure.  It is clear that $F(\Dirac)$ is still a weak
  equivalence.  Since both $F(\Dirac)$ and~$j$ are weak equivalences
  and $j\circ\Dirac'= F(\Dirac)$, it follows that~$\Dirac'$ is a weak
  equivalence.  This is a general fact about localization of
  triangulated categories.  Thus $\Dirac'\in\KK^G(\Proj,\C)$ is a
  Dirac morphism for~$G$.  Now let $\eta'\defeq F(\eta)\circ j$.  Then
  $\eta'\circ\Dirac'= F(\eta\Dirac)=\mathrm{id}_\Proj$ by
  construction.  Thus~$\eta'$ is a dual Dirac morphism.
\end{proof}

Theorem~\ref{uniformembeddingimpliesexistenceofgammaelement} now
follows by combining the above results.  If~$G$ uniformly embeds in a
Hilbert space, then it admits an a\nbd{}T\brd{}menable action on a
second countable compact space~$X$ by Theorem~\ref{embeddable_acts}.
Lemma~\ref{inducedactionisamenable} allows us to assume that~$X$ is
$H$\nbd{}equivariantly contractible for any compact subgroup
$H\subseteq G$.  The transformation group $G\cross X$ has a dual Dirac
morphism by Theorem~\ref{tustheorem}.  Hence so has~$G$ by
Proposition~\ref{pro:gamma}.

\end{document}